\documentclass[a4paper]{article}
\usepackage[english]{babel}
\usepackage{amsthm}  
\usepackage{amsmath} 
\usepackage{amssymb} 
\usepackage{mathabx} 
\usepackage{mathtools}
\usepackage{graphicx}
\usepackage{fancyhdr}
\usepackage{setspace}
\usepackage{lastpage}

\usepackage{subcaption}
\usepackage{adjustbox}
\usepackage{makecell}
\usepackage[numbers,sort&compress]{natbib}	

\usepackage{booktabs}
\usepackage{multirow}

\newcommand{\mytitle}{Global Deterministic Optimization with Artificial Neural Networks Embedded}
\newcommand{\myshorttitle}{Global Deterministic Optimization with ANNs Embedded}
\newcommand{\myauthor}{Artur M. Schweidtmann and Alexander Mitsos} 
\newcommand{\myauthorshort}{Schweidtmann and Mitsos}
\author{\myauthor}
\usepackage[colorlinks,linkcolor=blue,citecolor=blue,urlcolor=blue,pdftitle={\mytitle},pdfauthor={Artur M. Schweidtmann}]{hyperref}

{
	\fancyhead[L]{\myshorttitle}
	\fancyhead[R]{\small{\textit{\the\day.\the\month.\the\year}}}
	\fancyfoot[L]{\copyright \small{\textit{\myauthorshort}}}
	\fancyfoot[R]{Page \thepage\ of \pageref*{LastPage}}
	\cfoot{\textit{Pre-print}}
}

\fancypagestyle{firststyle}
{
	\fancyhead[L]{\copyright \small{\textit{\myauthorshort}}}
	\fancyhead[C]{\textit{Pre-print}}
	\fancyhead[R]{Page \thepage\ of \pageref*{LastPage}}
	\fancyfoot[C]{}
	\fancyfoot[R]{}
}

\addto\captionsenglish{%
}

\DeclareMathOperator{\sech}{sech}

\newtheoremstyle{myplain}
{} 
{} 
{\itshape} 
{} 
{} 
{} 
{ } 
{\textbf{\thmname{#1}\thmnumber{ #2}}\thmnote{ (#3)}} 

\newtheoremstyle{mydefinition}
{} 
{} 
{} 
{} 
{} 
{} 
{ } 
{\textbf{\thmname{#1}\thmnumber{ #2}}\thmnote{ (\textit{#3})}} 

\newtheoremstyle{myremark}
{} 
{} 
{} 
{} 
{} 
{} 
{ } 
{\textit{\thmname{#1}\thmnumber{ #2}\thmnote{ (#3)}}} 

\newtheoremstyle{mynote}
{} 
{} 
{} 
{} 
{} 
{:} 
{ } 
{\textit{\thmname{#1}\thmnumber{ #2}\thmnote{ (#3)}}} 

\newtheorem{proposition}{Proposition}[section]

\makeatletter
\let\@addpunct\@gobble
\g@addto@macro{\thm@space@setup}{\thm@headpunct{}} 
\makeatother

\begin{document}
	
	\thispagestyle{firststyle}
	\begin{flushleft}\begin{large}\textbf{\mytitle}\end{large} \end{flushleft}
	\myauthor
	
	\begin{flushleft}\begin{small}
			RWTH Aachen University\\
			AVT - Aachener Verfahrenstechnik\\
			Process Systems Engineering\\
			Forckenbeckstr. 51 \\
			D - 52074 Aachen
		\end{small}
	\end{flushleft}
	This is a pre-print of an article published in Journal of \textit{Optimization Theory and Applications}. The final authenticated version is available online at: \url{https://doi.org/10.1007/s10957-018-1396-0}.
	\newpage
\section{Abstract}
		Artificial neural networks (ANNs) are used in various applications for data-driven black-box modeling and subsequent optimization. Herein, we present an efficient method for deterministic global optimization of ANN embedded optimization problems. The proposed method is based on relaxations of algorithms using McCormick relaxations in a reduced-space [\textit{SIOPT}, 20 (2009), pp. 573-601] including the convex and concave envelopes of the nonlinear activation function of ANNs. The optimization problem is solved using our in-house global deterministic solver MAiNGO. The performance of the proposed method is shown in four optimization examples: an illustrative function, a fermentation process, a compressor plant and a chemical process optimization. The results show that computational solution time is favorable compared to the global general-purpose optimization solver BARON. \\[12pt]
		\textbf{Keywords:} Surrogate-based optimization, Multilayer perceptron, McCormick relaxations, Machine learning, Deep learning, Fermentation process, Compressor plant, Cumene process  \\[12pt]
		\textbf{AMS subject classification:} 90C26, 90C30, 90C90, 68T01
\section{Introduction}
\label{sec:Introduction}
A feed-forward ANN with one hidden layer can approximate any smooth function on a compact domain to an arbitrary degree of accuracy given a sufficient number of neurons \cite{Hornik.1989}. Thus, artificial neural networks (ANNs) are said to be universal approximators. Due to their excellent approximation properties, ANNs have attracted widespread interest in  various fields such as chemistry \cite{Gasteiger.1993}, chemical engineering \cite{AzlanHussain.1999}, pharmaceutical research \cite{AgatonovicKustrin.2000} and  biosorption processes \cite{WitekKrowiak.2014} for the black-box modeling of complex systems. Furthermore, various industry applications of ANNs exist, e.g., modeling and identification, optimization and process control \cite{Meireles.2003}. \newline \indent
In many of these applications, a model including ANNs is developed and then used for optimization. In biochemical engineering for example, ANNs are commonly used to model fermentation processes \cite{DelRioChanona.2017} and optimize their output to identify promising operating points for further experiments \cite{Cheema.2002,Desai.2008,Nagata.2003}. ANNs are also used as surrogate models for the optimization of chemical processes \cite{Fahmi.2012,Nascimento.2000,Nascimento.1998,Chambers.2002,Henao.2010,Henao.2011,SantAnna.2017,Smith.2013,Henao.2012,Kajero.2017,Lewandowski.1998,GutierrezAntonio.2016,Chen.2002,Dornier.1995}. \newline \indent 
Most of the previous optimizations with ANNs embedded rely on local optimization techniques. Fernandes, for instance, optimized the product concentrations of a Fischer-Tropsch synthesis using the quasi-Newton method and finite-differences \cite{Fernandes.2006}. Other authors used local solves like DICOPT \cite{Grossmann.2002} (e.g., \cite{Fahmi.2012,Henao.2010,Henao.2011,Henao.2012}), sequential quadratic programming (e.g., \cite{SantAnna.2017}) or CONOPT \cite{Drud.1994} (e.g., \cite{Henao.2012}). These local methods can yield suboptimal solutions when the learned input-output relation is multi-modal. In addition, the activations functions involved in ANNs are usually nonconvex. A few researchers have addressed this problem by using stochastic global methods such as genetic algorithms (e.g., \cite{Nandi.2001,Cheema.2002,Desai.2008,Nagata.2003,GutierrezAntonio.2016,Chen.2002}) or brute-force grid search (e.g., \cite{Nascimento.1998,Nascimento.2000,Chambers.2002}). These methods cannot guarantee global optimality. Global deterministic optimization of an ANN embedded problem was done by Smith \textit{et al.} (2013) \cite{Smith.2013} who used BARON \cite{Ryoo.1996,Tawarmalani.2004,Tawarmalani.2005} to optimize an ANN with one hidden layer and three neurons that emulates a flooded bed algae bioreactor. As the hyperbolic tangent activation function is currently not available in BARON, Smith \textit{et al.} reformulated the activation function using exponential functions. Henao argues in his Ph.D. thesis \cite{Henao.2012} that state-of-the-art global optimization solvers such as BARON are at that time, i.e., 2012 limited to small to medium sized problems and he suggests a piecewise linear approximation of the hyperbolic tangent function. The literature shows that there remains a need for an efficient deterministic global optimization algorithm of problems with ANNs embedded. \newline \indent
There have been many other efforts to combine surrogate modeling with (global) optimization. For instance, Gaussian processes (GPs) have been used in the field of Bayesian optimization for optimization of expensive-to-evaluate black-box functions (e.g., \cite{Bradford.2017,Forrester.2009,Jones.2001,Sacks.1989,Shahriari.2016}). In the adaptive approach, GPs emulate black-box objectives, e.g., measurements of a chemical experiment \cite{Schweidtmann.2018}, in every iteration to assist a following query point selection. ALAMO is another adaptive sampling approach \cite{Wilson.2017,Cozad.2015,Cozad.2014} that aims the development of simple surrogate models in the light of small data sets.
This contribution focuses on the development of a global deterministic optimization approach for problems including (given) well-established ANN surrogate models. In contrast to standard GPs, ANNs can have multiple outputs and can handle large training sets (i.e., GP training requires N by N matrix inversions where N is the number of training points). \newline \indent
In deterministic global optimization, convex relaxations are derived and relaxed problems are solved using branch-and-bound (B\&B) \cite{Falk.1969,Horst.1996} or related algorithms \cite{Misener.2014,Ryoo.1996,Tawarmalani.2004,Tawarmalani.2005}. 
In most state-of-the-art deterministic global optimization solvers such as BARON \cite{Ryoo.1996,Tawarmalani.2004,Tawarmalani.2005}, ANTIGONE \cite{Misener.2014} and SCIP \cite{Maher.2017} the complete set of constraints and variables is handed to the optimization algorithm that builds convex relaxations. This method is referred to as the full-space (FS) formulation. Another possibility is the reduced-space (RS) formulation. The key idea herein is that the optimizer does not see all variables and constraints \cite{Epperly.1997,Mitsos.2009,Scott.2011,Bongartz.2017,Huster.2017}. \newline \indent
For the deterministic global optimization of ANNs, the FS inherently results in a large-scale optimization problem because the ANN network structure includes multiple (hidden) layers, neurons and network equations. Additionally, auxiliary variables used in state-of-the-art solvers for the relaxations increase the problem size. Globally solving large-scale optimization problems is difficult because of the exponential worst-case runtime of the B\&B algorithm. Moreover, in standard solvers the user has to provide tight variable bounds which are difficult to provide for some network variables. A possible disadvantage of the RS formulation is that the propagation of relaxations through large ANNs may result in weak relaxations. \newline \indent
One possibility to construct relaxations in the RS setup are McCormick relaxations \cite{McCormick.1976,McCormick.1983} that have shown favorable convergence properties \cite{Bompadre.2012,Najman.2016} and have been extended to the relaxation of multivariate functions \cite{Tsoukalas.2014} and bounded functions with discontinuities \cite{Wechsung.2014}. Recently, a differentiable modification of the relaxations has been proposed as well \cite{Khan.2017}. The global optimization using McCormick relaxations has been applied to several problems such as flowsheet optimization \cite{Bongartz.2017,Bongartz.2017b,Huster.2017} and the solution of ODE and nonlinear equation systems \cite{Wechsung.2015,Scott.2011,Stuber.2015}. Besides the direct utilization of McCormick relaxations in the original variables, the method was also used in the development of the auxiliary variable method (AVM) that is applied in state-of-the-art solvers \cite{Misener.2014,Smith.1997,Tawarmalani.2002}. \newline \indent
In this contribution, ANNs are optimized in a RS which significantly reduces the dimensionality of the optimization problems. More specifically, the equations that describe the ANN and the corresponding variables are hidden from the B\&B solver.
In order to construct tight convex and concave relaxations of the ANNs, we utilize the convex and concave envelopes of the activation function (shown in Appendix \ref{sec:Convex and concave envelopes of the hyperbolic tangent activation function}) and automatic construction of McCormick relaxations \cite{Mitsos.2009}. \newline \indent 
In the remainder of this paper, first an overview about multilayer perceptron ANNs is provided (Section \ref{sec:Background}). In Section \ref{sec:Relaxationsofneuralnetworks}, the optimization problem formulation and the relaxation of ANNs are proposed. Further, implementation details are provided. In Section \ref{sec:Numericalresults}, the proposed method is applied to four numerical examples illustrating its potential and compared to a state-of-the-art general purpose optimization algorithm BARON in terms of computational (CPU) time. In Section \ref{sec:Conclusionandfuturework}, advantages, limitations and prospective utility of the proposed method are discussed.
\section{Background on Multilayer Perceptrons}
\label{sec:Background}
\begin{figure}[ht]
	\begin{subfigure}[c]{0.5\textwidth}
		\includegraphics[width=2.0 in]{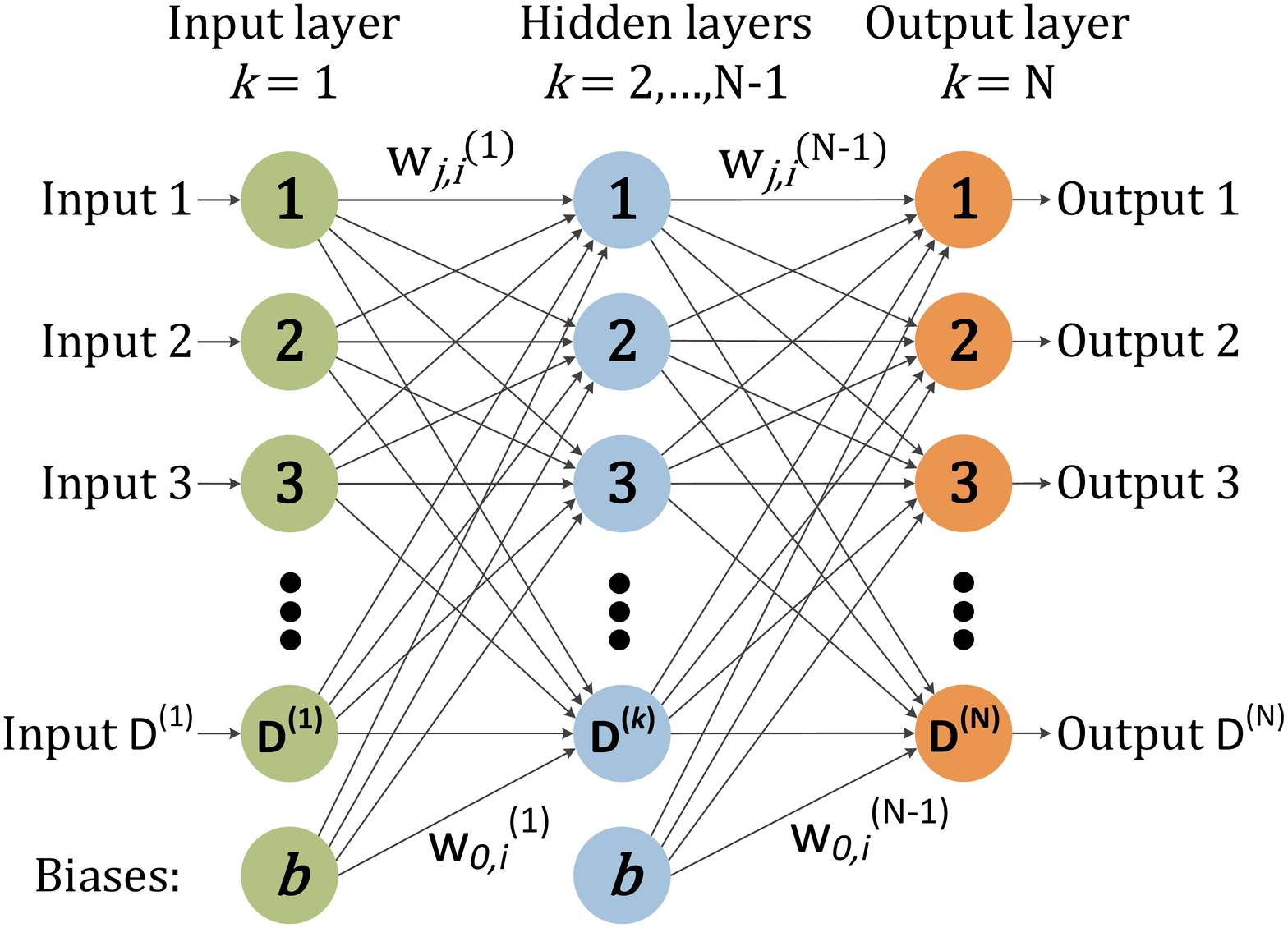}
		\subcaption{The complete network}
	\end{subfigure}
	\begin{subfigure}[c]{0.5\textwidth}
		\includegraphics[width=2.0 in]{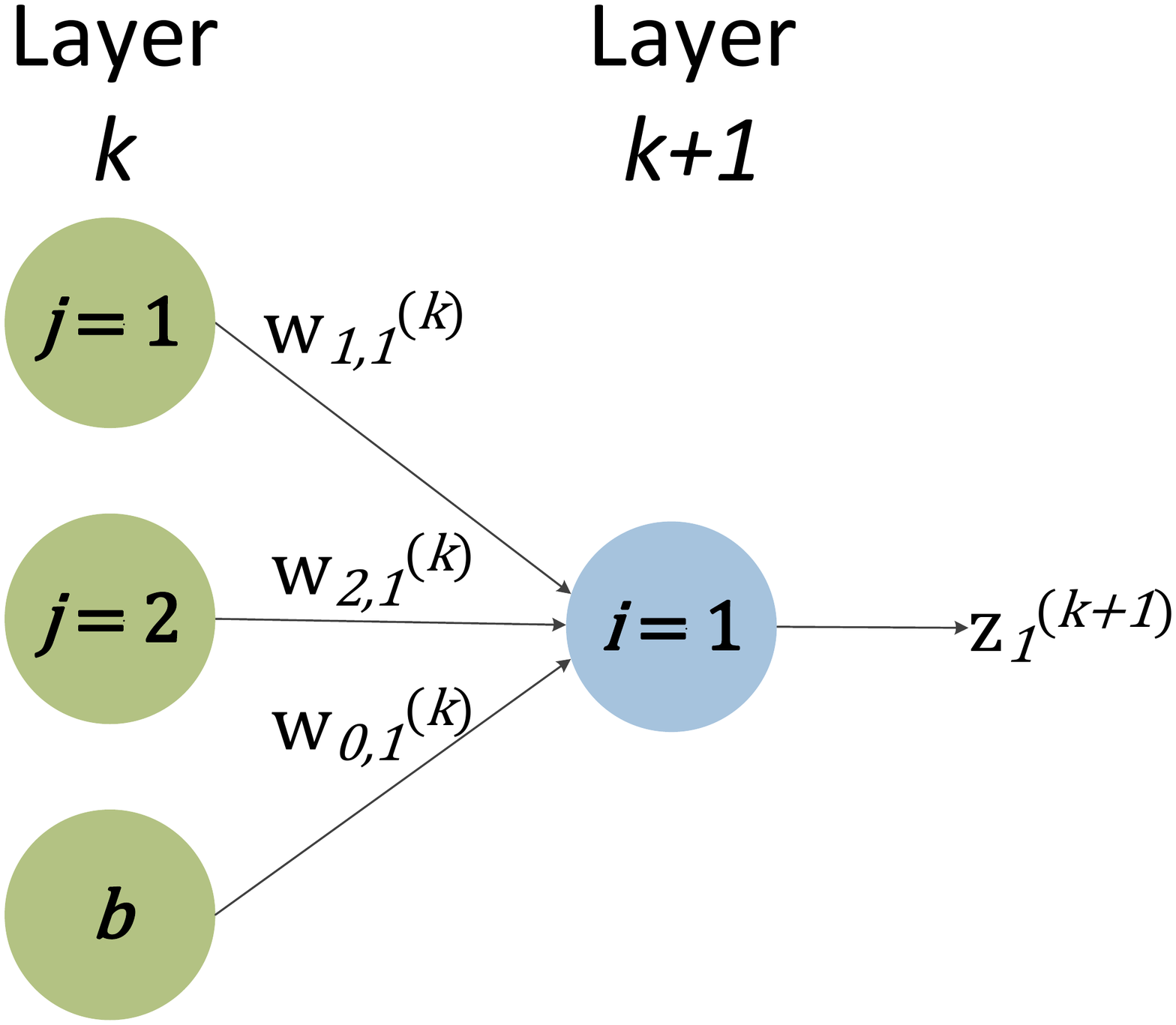}
		\subcaption{A detail of the network}
	\end{subfigure}
	\caption{Graphical illustration of a multilayer perceptron artificial neural network as a directed acyclic graph}
	\label{fig:1}  
\end{figure}

In this subsection, the multilayer perceptron (MLP), also known as feed-forward ANN, is briefly introduced. More detailed information about MLPs can be found in the literature (e.g., \cite{Bishop.2009}). As depicted in Figure \ref{fig:1}, the MLP can be illustrated as a directed acyclic graph connecting multiple layers, $k$, of neurons, which are the nodes of the graph. It consists of an input layer ($k=1$), a number of hidden layers ($k = 2,..,N-1$) and an output layer ($k = N$). Each connection between a neuron $j$ of layer $k$ and neuron $i$ of layer $k+1$ is associated with a weight $w^{(k)}_{j,i}$. The output, $z^{(k+1)}_{i}$, of a neuron $i$ in layer $k$ is given by 
\begin{equation}\label{eqn:ANN}
z^{(k+1)}_{i} = h_{k+1} \left( \sum_{j = 1}^{D^{(k)}} (w^{(k)}_{j,i} z^{(k)}_{j}) + w^{(k)}_{0,i} \right)
\end{equation}
where $h_{k}:\mathbb{R} \to \mathbb{R}$ is the activation function, $D^{(k)}$ is the number of neurons in each layer and $w^{(k)}_{0,i}$ is a bias parameter. Thus, the argument of $h$ is a linear combination of the outputs of the previous layer. The output value of neurons in the input layer ($z^{(1)}_{i}$) correspond to the network input variables. Following \eqref{eqn:ANN} for all neurons and layers, the outputs (also called predictions) of the network ($z^{(N)}_{i}$) can be computed. 
 \newline \indent
A commonly-used activation function of the hidden layers is the hyperbolic tangent function ($h_k(x)=\tanh(x)$). In the output layer, the identity activation function ($h_N(x)=x$) is typically used for regression problems. Thus, this work focuses on the hyperbolic tangent function. However, the proposed method can be easily applied to other common activation functions like the sigmoid function (see Appendix \ref{sec:Convex and concave envelopes of the sigmoid activation function}).
\section{Method}
\label{sec:Relaxationsofneuralnetworks}
In this section, the optimization problem formulations and the McCormick relaxation of MLPs are proposed. Further, details about the implementation are provided.
\subsection{Optimization Problem Formulations}
Optimization problems which use MLPs as surrogate models often have a particular structure. Usually, one or more of the input variables, $\textbf{x}$, correspond to degrees of freedom of the problem and can be chosen within given bounds $D = [\textbf{x}^{L},\textbf{x}^{U}]$. Once the input variables are fixed, the dependent variables in the networks, $\textbf{z}$, can be determined by solving the nonlinear network equations $\textbf{h}(\textbf{x},\textbf{z})=\textbf{0}, \quad \textbf{h}:D \times \mathbb{R}^{n_z} \to \mathbb{R}^{n_z}$. The nonlinear network equations of MLPs can also be formulated explicitly in the outputs as shown in Subsection \ref{subsec:Reduced-space formulation}. The outputs $\textbf{y}$ of the networks can be understood as network variables of the output layers. Thus, the output variables are a subset of the network variables with $n_y < n_z$. The objective function, $f(\textbf{x},\textbf{y}), \quad f:D \times \mathbb{R}^{n_y} \to \mathbb{R}$, usually depends on the networks inputs and outputs but not on the remaining network variables. Similarly, the feasible region is often constrained by inequalities $\textbf{g}(\textbf{x},\textbf{y}) \le \textbf{0}, \quad \textbf{g}: D \times \mathbb{R}^{n_y} \to \mathbb{R}^{n_g}$ which also depend only on the network inputs and outputs.
\subsubsection{Full-space formulation}
\label{subsec:Full-space formulation}
In the FS formulation, the input and the network variables are optimization variables and the problem can be formulated as follows:
\begin{align}
& \underset{\textbf{x} \in D, \textbf{z} \in Z}{\min} 	& f(\textbf{x},\textbf{y}) 			\notag		\\
& \text{s.t.}  											& \textbf{h}(\textbf{x},\textbf{z})	& = \textbf{0} 	\tag*{(FS)}	\\
& 														& \textbf{g}(\textbf{x},\textbf{y}) 	& \le \textbf{0} \notag		
\end{align}
In this case, a global B\&B solver requires bounds on $\textbf{x}$ and $\textbf{z}$. As mentioned in the introduction, the FS formulation is utilized by common general-purpose deterministic global optimization solvers.
\subsubsection{Reduced-space formulation}
\label{subsec:Reduced-space formulation}
In MLPs the network equations $\textbf{h}(\textbf{x},\textbf{z})=\textbf{0}$ can be reformulated as an explicit function $\hat{\textbf{y}}: D \to \mathbb{R}^{n_y}$ with $\textbf{y} = \hat{\textbf{y}}(\textbf{x})$.
Thus, the dependent network variables can be eliminated from the FS optimization formulation (c.f. \cite{Bongartz.2017}). The resulting optimization problem can be formulated as follows:
\begin{align}
& \underset{\textbf{x} \in D}{\min} 					& f(\textbf{x},\hat{\textbf{y}}(\textbf{x}))		 			\tag*{(RS)}		\\
& \text{s.t.}  											& \textbf{g}(\textbf{x},\hat{\textbf{y}}(\textbf{x}))		& \le \textbf{0}  	\notag
\end{align}
The RS formulation operates only in the domain $D$ of the optimization variables $x$. Further, the equality constraints which describe the network equations are not visible to the optimization algorithm anymore. However, it is important to note that the complete elimination of equality constraints is not necessarily possible when optimizing arrangements of several MLPs or hybrid modeling formulations. In these cases, some equality constraints and additional variables might remain in the RS formulation (e.g., for connecting different MLPs with a recycle) similar to, e.g., \cite{Bongartz.2017}. As an option, these can also be relaxed using extensions of the McCormick approach to implicit functions \cite{Stuber.2015,Wechsung.2015}.
\subsection{Relaxations of Artificial Neural Networks}
In order to solve the RS optimization problem using the B\&B algorithm, lower-bounds of the MLPs have to be derived. In this paper, the lower bounds are calculated by automatic propagation of McCormick relaxations and natural interval extensions \cite{Moore.1979,Hofschuster.1998} using the open-source software MC\texttt{++} \cite{Chachuat.2014,Chachuat.2015} (see Section \ref{sec:Implementation}). As described in Section \ref{sec:Background}, the outputs of MLPs are compositions of the activation functions and affine functions that connect the neurons. Thus, the only nonlinearities of MLPs are its activation functions. When using MLPs as a surrogate model, the hyperbolic tangent function ($\tanh$) is commonly used as the activation function in the hidden layers and the identity is usually used as the activation function in the output layer \cite{Bishop.2009}. Thus, the tightness of MLP relaxations is considerably influenced by the relaxations of the hyperbolic tangent function. \newline \indent
In general, McCormick relaxations do not provide the  envelopes (tightest possible relaxations). The envelopes of the hyperbolic tangent function are derived in Appendix \ref{sec:Convex and concave envelopes of the hyperbolic tangent activation function}. They are once continuously differentiable (Proposition \ref{prop:smoothness of hyperbolic tangent relaxations}) and strictly monotonically increasing (Proposition \ref{prop:monotonicity of hyperbolic tangent relaxations}). \newline \indent
In the proposed framework, the envelopes of the hyperbolic tangent activation function are propagated through the network equations to derive concave and convex relaxations of the outputs of MLPs. The resulting relaxations of MLPs are once continuously differentiable as shown in the following proposition.
\begin{proposition} \label{prop:smoothness of MLP relaxations}
	(smoothness of MLP relaxations). Consider an MLP using (exclusively) the hyperbolic and identity activation functions. The McCormick relaxations built with the envelopes of the hyperbolic tangent function (Appendix \ref{sec:Convex and concave envelopes of the hyperbolic tangent activation function}) are once continuously differentiable.
\end{proposition}
{\it Proof}
	Consider the output of a neuron in a hidden layer as a function of its inputs. As the convex and concave relaxations of the hyperbolic tangent activation function are monotonically increasing (see Proposition \ref{prop:monotonicity of hyperbolic tangent relaxations}), Corollary 3 by \cite{Tsoukalas.2014} holds. Thus, the convex and concave relaxation of the output of the neuron are compositions of the corresponding relaxation of the hyperbolic tangent function and the corresponding relaxations of the inputs to the neuron. As the envelopes of the hyperbolic tangent function are once continuous differentiable (see Proposition \ref{prop:smoothness of hyperbolic tangent relaxations}), the relaxations of the MLPs are once continuous differentiable by chain rule. 
\qed
It should be noted that the relaxations are not $C^2$ because the envelopes of the hyperbolic tangent function are not $C^2$ (see Proposition \ref{prop:smoothness of hyperbolic tangent relaxations}). However, Proposition \ref{prop:lipschitz} shows that the first derivative of the convex and concave envelopes of the hyperbolic tangent function, $\frac{ \text{d}(F^{cv})}{\text{d} x}$ and $\frac{ \text{d}(F^{cc})}{\text{d} x}$, are Lipschitz continuous (sometimes referred to as $C^{1,1}$). Recall that we use the hyperbolic tangent function because it is the most commonly used and similar results can be derived for other activition functions. \newline \indent
Proposition \ref{prop:smoothness of MLP relaxations} has a beneficial consequence for the global optimization of MLPs. In general, McCormick relaxations are continuous but not differentiable ($C^0$). This necessitates nonsmooth algorithms, e.g., bundle methods \cite{Bertsekas.2003,Bertsekas.2015} or linearization-based methods \cite{Mitsos.2009}. However, when optimizing MLPs, a $C^{1,1}$ optimization algorithm can be used which is potentially more efficient. \newline \indent
In the general-purpose optimization solvers BARON, ANTIGONE and SCIP, the hyperbolic tangent function is currently not available and interfaces do not allow the user to define functions and/or relaxations. In order to compare the performance of the proposed method to these algorithms, the hyperbolic tangent activation function can be reformulated. As the tightness of the McCormick relaxations of different reformulations can differ, the relaxations of four common reformulations $F_1, F_2, F_3, F_4: \mathbb{R} \to \mathbb{R}$ with $F_1(x) \coloneqq \frac{e^x - e^{-x}}{e^x + e^{-x}},~	F_2(x)  \coloneqq \frac{e^{2x} - 1}{e^{2x} + 1},~ F_3(x) \coloneqq 1 - \frac{2}{e^{2x} + 1},~ F_4(x) \coloneqq \frac{1 - e^{-2x}}{1 + e^{-2x}}$ are compared in Appendix \ref{sec:McCormick relaxations of reformulations of the hyperbolic tangent activation function}. The exemplary plot of the McCormick relaxations in the appendix shows that the convex and concave relaxations of $F_1$, $F_2$ and $F_4$ are weaker than the ones of $F_3$. Further, the convex and concave relaxations of $F_1$, $F_2$ and $F_4$ are not differentiable. Finally, it can be observed that relaxations of the reformulations are considerably weaker than the envelopes.
\subsection{Implementation}
\label{sec:Implementation}
The proposed work uses the in-house \textbf{M}cCormick based \textbf{A}lgorithm for mixed \textbf{i}nteger \textbf{N}onlinear \textbf{G}lobal \textbf{O}ptimization (MAiNGO), a B\&B optimization solver in C\texttt{++} \cite{Bongartz.2018}. Herein, a best-first heuristic and bisection along the longest edge is used for branching. For lower bounding, a convex underestimation of the problem is found by automatic propagation of McCormick relaxations using the open-source software MC\texttt{++} v2.0 \cite{Chachuat.2014,Chachuat.2015}. The necessary interval extensions are provided by FILIB\texttt{++} v3.0.2 \cite{Hofschuster.1998}.  The convex relaxations of the constraints and the objective are linearized with the use of subgradients at the centerpoint of each node and the resulting linear program (LP) is solved by CPLEX v12.5 \cite{InternationalBusinessMachies.2009}. Further, optimization-based bound tightening is used that is improved by filtering bounds technique with a factor of 0.1 as described in \cite{Gleixner.2017} and bound tightening based on the dual multipliers returned by CPLEX is used \cite{Ryoo.1995,Locatelli.2013}. For upper bounding, the problem is locally optimized using the SLSQP algorithm \cite{Kraft.1988,Kraft.1994} in the NLopt library v2.4.2 \cite{Johnson.2016}. The necessary derivatives are provided by the FADBAD\texttt{++} tool for automatic differentiation \cite{Bendtsen.2012}. Furthermore, recently developed heuristics for tighter McCormick relaxations are used \cite{Najman.2017b}. \newline \indent 
The convex and concave envelopes of $\tanh$ are added to the MC\texttt{++} library. In order to compute those envelopes, the equations \eqref{eqn:x_c^u} and \eqref{eqn:x_c^o} have to be solved numerically (see Appendix \ref{sec:Convex and concave envelopes of the hyperbolic tangent activation function}). For this purpose, the Newton method is used and analytical gradient information is supplied. \newline \indent
The ANNs in this work are fitted in the Neural Network Toolbox in MATLAB. 
\section{Numerical Results}
\label{sec:Numericalresults}
In this section, the numerical results of four case studies are presented. The performance of the proposed method is compared to the general-purpose optimization solver BARON 17.4.1. using GAMS 24.8.5. All numerical examples were run on one thread of an Intel Xeon CPU E5-2630 v2 with 2.6 GHz, 128 GB RAM and Windows Server 2008 operating system.
\subsection{Illustrative Example \& Scaling of the Algorithm}
\label{subsec:Numerical Example 1}
\begin{figure}[ht]
	\begin{subfigure}[c]{0.5\textwidth}
		\includegraphics[width=2.0 in]{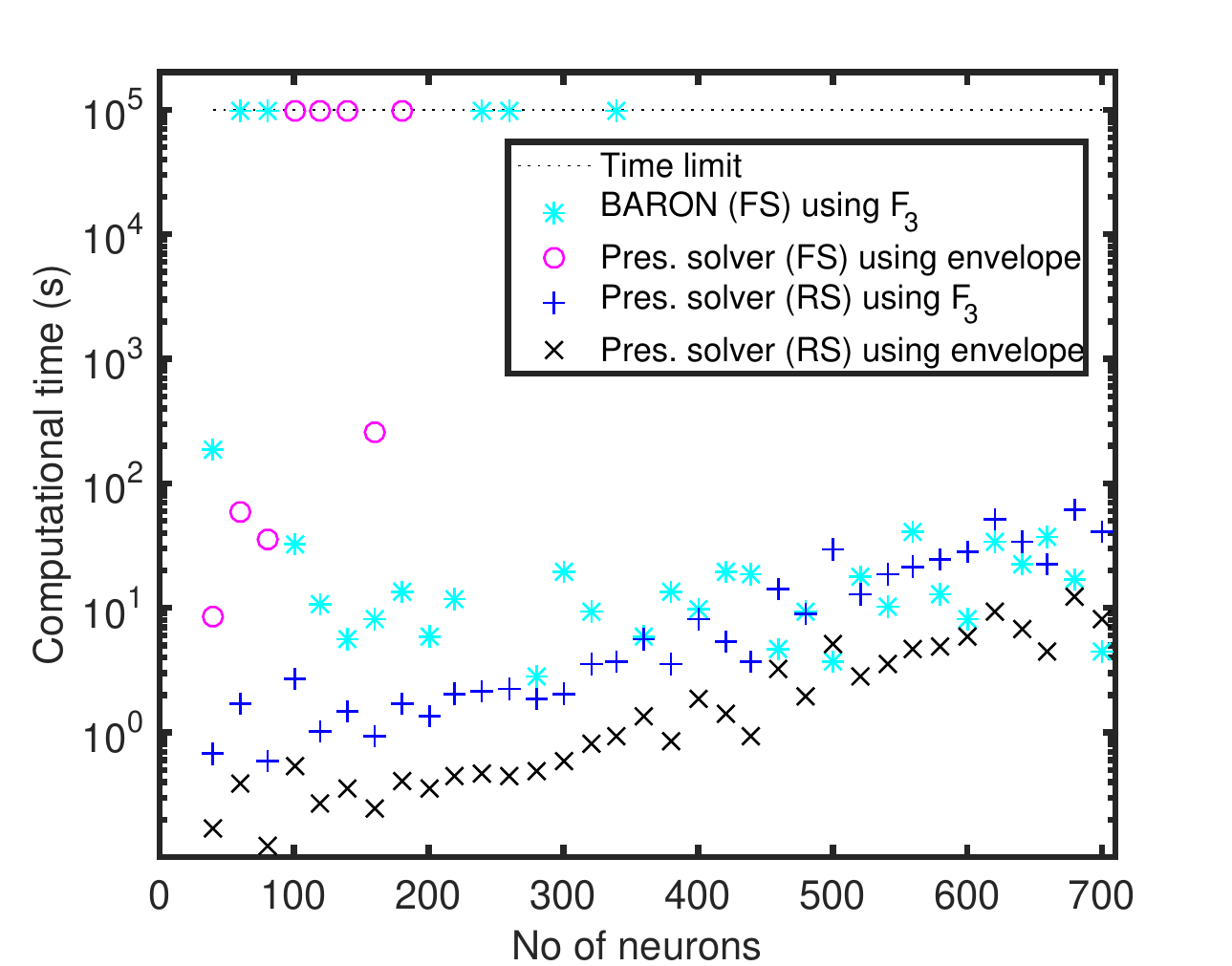}
		\subcaption{MLP with one hidden layer and varying number of neurons}
	\end{subfigure}
	\begin{subfigure}[c]{0.5\textwidth}
		\includegraphics[width=2.0 in]{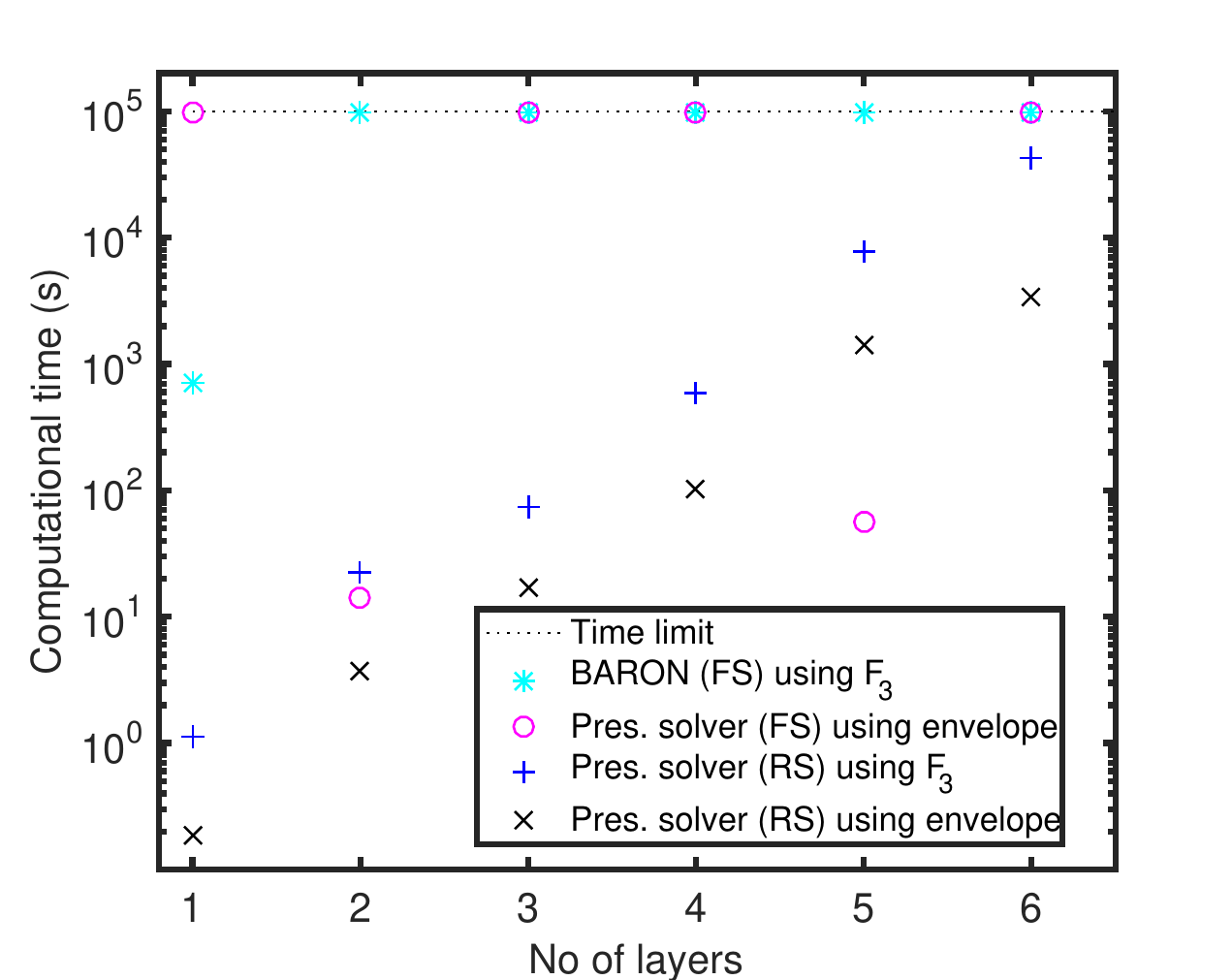}
		\subcaption{MLP with 20 neurons in each layer and varying number of hidden layers}
	\end{subfigure}
	\caption{Graphical illustration of the computational time over the ANN size for the peaks function optimization problem (Subsection \ref{subsec:Numerical Example 1})}
	\label{fig:2}
\end{figure}

In the first illustrative example, a two-dimensional mathematical test function is learned by a MLP that is subsequently optimized. The peaks function is provided by Matlab (\textit{peaks()}) and is given by $f_{peaks}: \mathbb{R}^2 \to \mathbb{R}$ with 
\begin{equation}
f_{peaks}(x_1,x_2) =  3 (1-x_1)^2 \cdot e^{-x_1^2 - (x_2+1)^2} - 10 \cdot (\frac{x_1}{5} - x_1^3 - x_2^5) \cdot e^{-x_1^2-x_2^2}- \frac{e^{-(x_1+1)^2 - x_2^2} }{3} 
\end{equation} \indent 
The function has multiple local optima on $D = \{x_1, x_2 |-3 \leq x_1, x_2 \leq 3 \}$. The known unique global minimizer of the test function is $\min_{x_1,x_2 \in D} f_{peaks}(x_1,x_2) = -6.551$ at $(x_1^*,x_2^*) = (0.228,-1.626)$. \newline \indent 
To learn the peaks function, a Latin hypercube (LHC) sampling technique is used to generate a set of 500 points on $D$. As mentioned in the introduction, a MLP with one hidden layer is a universal approximator. Thus, a network architecture with one hidden layer is chosen in this example. The MLP constitutes two neurons in the input layer and one neuron in the output layer. In order to avoid over- and under-fitting, $D$ is randomly divided into a training, a validation and a test sets with the respective size ratios of $0.7:0.15:0.15$. The weights of the network are fitted in the Matlab Neural Network Toolbox by minimizing the mean squared error (MSE) on the training set using Levenberg-Marquardt algorithm and the early stopping procedure \cite{Bishop.2009}. To obtain a suitable number of neurons in the hidden layer the training is repeated using different number of neurons. The configuration with the lowest MSE on the test set determines the number of neurons to be 47. The training results in a MSE of $6.8 \cdot 10^{-5}$, $3.9 \cdot 10^{-4}$  and $4.6 \cdot 10^{-4}$  on the training, validation and test set respectively. \newline \indent 
After training, the MLP is optimized using the proposed methods. The absolute optimization tolerance is set to $\epsilon_\text{tol} = 10^{-4}$ because more accurate solution of the optimization problem is not sensible due to the prediction error of the MLP. The relative tolerance is set to its minimum value ($10^{-12}$) and is thus not active. This way the optimization result is limited by the prediction accuracy only. Further, a time limit of 100,000 seconds (about 28 hours) is set for the optimization. \newline \indent
Table \ref{tab:NumericalResults_Peaks} summarizes the results of the optimization runs. The results show that the problem formulation has a significant influence on the problem size. In comparison to the FS formulations, the RS formulation reduces the number of the optimization variable from 101 to 2 and eliminates all 99 equality constraints. The RS formulation requires over 190 times less CPU time compared to the FS formulation when using the presented solver. Further, different relaxations of the activation function yield different solution times. The envelope of the activation function accelerates the optimization significantly compared to all reformulations. When the presented solver is used in the FS formulation, the weak relaxations of the reformulations even lead to a variable overflow error in the exponential function in FILIB\texttt{++}. As expected from the comparison of the relaxations (see Appendix \ref{sec:McCormick relaxations of reformulations of the hyperbolic tangent activation function}), reformulation $F_3$ outperforms the other reformulations in terms of computational efficiency. In general, all RS optimizations converged to the same optimal solution ($-6.563$). The solution is about 0.18\% different from the global optimizer of the underlying peaks function. This is within the expected accuracy of the MLP prediction.
\begin{table}[ht]
	\centering
	\caption{Numerical results of the peaks function optimization (Subsection \ref{subsec:Numerical Example 1})}
	\label{tab:NumericalResults_Peaks}
	\resizebox{\columnwidth}{!}{
		\begin{tabular}{lllrrr}
			\hline
			& Problem size & Solver                      & CPU time {[}s{]} & Iterations & Abs. gap\\ \hline
			\multirow{9}{*}{(FS)} 
			&  101 variables & BARON ($F_1$)                &     1,727.29              &  55,047  	 	&  $< \epsilon_\text{tol}$    \\
			& 99 equalities & BARON ($F_2$)                	&       387.44              &    13,329   	&  $< \epsilon_\text{tol}$   \\
			& 0 inequalities & BARON ($F_3$)                &        21.08              &   1,223   	&  $< \epsilon_\text{tol}$    \\
			&              	& BARON ($F_4$)                	&   100,000.00              &    2,732,488  & 0.0003     \\
			
			&  				& Presented solver ($F_1$)      &  	\multicolumn{3}{c}{--variable overflow--}   \\
			&             	& Presented solver ($F_2$)      &  	\multicolumn{3}{c}{--variable overflow--}    \\
			&             	& Presented solver ($F_3$)     	& 	\multicolumn{3}{c}{--variable overflow--}     \\
			&              	& Presented solver ($F_4$) 		&  	\multicolumn{3}{c}{--variable overflow--}   \\
			&              	& Presented solver (envelope) 	&   84.19               &  1     & $< \epsilon_\text{tol}$   \\
			
			\hline
			\multirow{5}{*}{(RS)}  
			& 2 variables & Presented solver ($F_1$)          &  18.61                & 8,825    &  $< \epsilon_\text{tol}$    \\
			& 0 equalities             & Presented solver ($F_2$)          &  23.06                & 17,785    &  $< \epsilon_\text{tol}$     \\
			& 0 inequalities            & Presented solver ($F_3$)     &  0.73         &  933    & $< \epsilon_\text{tol}$     \\
			&              & Presented solver ($F_4$) &   10.50               & 11,857     & $< \epsilon_\text{tol}$    \\
			&              & Presented solver (envelope) &   0.44               &  207     & $< \epsilon_\text{tol}$   \\
			\hline
			\multicolumn{6}{l}{variable overflow: Crashed due to variable overflow in FILIB\texttt{++}.} \\
		\end{tabular}
	}
\end{table}  \newline \indent 
In a second step, the peaks function is used to illustrate the scaling of the optimization algorithm with network size. Therefore, networks with different sizes are trained on a LHC set of 2,000 points. Subsequently, the networks are optimized. In Subfigure \ref{fig:2} (a), the number of neurons in the hidden layer of a shallow MLP is varied from 40 to 700 and the computational time for optimization is depicted.  The RS formulation using the envelope and the reformulation $F_3$ shows a consistent increase of CPU time with the number of neurons whereas the CPU time of BARON using $F_3$ does not show a clear trend. For instance, BARON did not converge within 100,000 seconds for networks with 60, 80, 240, 260  and 340 neurons but shows good performance for a network with 700 neurons. The optimization using the FS formulation and the presented solver shows the weakest performance. As most of these optimizations with more than 100 neurons did not converge within 100,000 seconds, this optimization was not executed for networks with more than 100 neurons. The presented solver using the RS formulation and the envelope shows consistently one of the best performances. It should be noted that the number of neurons is unreasonable high for this illustration which leads to overfitting. It can further be noted that the network sizes (up to 700 neurons) exceed the network size of the case study by Smith \textit{et al.} (2013) \cite{Smith.2013} (3 neurons) drastically.
 \newline \indent 
In Subfigure \ref{fig:2} (b), the number of neurons in each hidden layer is fixed to 20 and the number of hidden layers is varied from 1 to 6.  The results show that the CPU times of the presented solver using the RS formulation  increases approximately exponentially with the number of layers. Further, the utilization of the envelopes has always an advantage over the reformulation $F_3$. The solution of the FS formulation using either BARON or the presented solver does not converge to an optimal solution within the time limit for most cases. In general, it is apparent that deep networks require more CPU time for optimization than shallow networks with the same or even larger number of neurons.
\subsection{Fermentation Process}
\label{subsec:Numerical Example 2}
In the second subsection, a fermentation of glucose to gluconic acid is learned from experimental data and optimized. The example is based on and compared to the work of Cheema \textit{et al.}  \cite{Cheema.2002} where an MLP is learned from experimental data and optimized using stochastic optimization approaches. The numerical example is relevant because mechanistic models of fermentation processes are often not available and surrogate-based optimization can help to identify promising operating conditions. \newline \indent
The inputs of the MLP and degrees of freedom are the glucose concentration ($x_1$ in [g/L]), the biomass concentration ($x_2$ in [g/L]) and the dissolved oxygen ($x_3$ in [mg/L]), whereas the output of the MLP is the concentration of gloconic acid ($y$ in [g/L]). The objective of the optimization problem is to maximize the yield of gluconic acid defined as $y_{gl} = \frac{100 y}{1.088 x_1}$. As the network weights used by Cheema \textit{et al.} \cite{Cheema.2002} are not available, we trained a MLP using the same structure, training and validation data, and training algorithm. \newline \indent 
After the training, we performed deterministic global optimization to maximize the gluconic acid yield. The result of the global deterministic optimization is $x_1 = 156.466$ g/L, $x_2 = 3$ g/L, $x_3 =  57.086 $ mg/L, $y = 170.127$ g/L and $y_{gl} = 99.937$ which is similar to the literature result. However, it is important to note that the underlying network is rebuild based on the data and it is very unlikely that the network is identical to the one from literature.  The main difference to the results from literature is that all deterministic methods in this work converge to the same solution whereas the stochastic literature method shows considerable variations in their solutions, i.e., even the top three out of hundreds of executions of the genetic algorithm (GA) and the simultaneous perturbation stochastic approximation (SPSA) show variations \cite{Cheema.2002}. \newline \indent 
Due to the small network size, the problem size of the RS is only slightly smaller than the FS formulation (compare Table \ref{tab:NumericalResults_Fermentation}). The CPU time of the presented solver using the envelope and RS formulation is about 5.4 times faster than the FS formulation. Thus, the RS formulation is favorable with a CPU time of 0.11 seconds. Again, different reformulation of the activation function yield different CPU times. In particular, the CPU time of reformulations $F_1$ and $F_2$ in the RS are a unexpectedly large. When considering only reformulations $F_1$ and $F_2$, BARON converges significantly faster to an optimal solution. Further, the CPU time of BARON using different reformulations does not seem to directly relate to the tightness of their McCormick relaxations.
\begin{table}[ht]
	\centering
	\caption{Numerical results of the fermentation process optimization (Subsection \ref{subsec:Numerical Example 2})}
	\label{tab:NumericalResults_Fermentation}
	\resizebox{\columnwidth}{!}{
		\begin{tabular}{lllrrr}
			\hline
			& Problem size & Solver                      & CPU time {[}s{]} & Iterations & Abs. gap \\ \hline
			\multirow{9}{*}{(FS)} 
			&  13 variables  	& BARON ($F_1$)                &  14.38                &  8,233 &    $< \epsilon_\text{tol}$     \\
			& 10 equalities 	& BARON ($F_2$)                &    0.44               &    5    &  $< \epsilon_\text{tol}$  \\
			& 0 inequalities    & BARON ($F_3$)                &    0.56               &   7   &   $< \epsilon_\text{tol}$   \\
			&              		& BARON ($F_4$)                &    0.48               &    15  &     $< \epsilon_\text{tol}$  \\
			
			&  				& Presented solver ($F_1$)      &  	\multicolumn{3}{c}{--variable overflow--}   \\
			&             	& Presented solver ($F_2$)      &  	\multicolumn{3}{c}{--variable overflow--}    \\
			&             	& Presented solver ($F_3$)     	& 	\multicolumn{3}{c}{--variable overflow--}     \\
			&              	& Presented solver ($F_4$) 		&  	\multicolumn{3}{c}{--variable overflow--}   \\
			&              	& Presented solver (envelope) 	&   0.59            &   157     & $< \epsilon_\text{tol}$   \\
			\hline
			\multirow{5}{*}{(RS)}  
			& 3 variables & Presented solver ($F_1$)          &  1,100.43                & 60,752,735   &   $< \epsilon_\text{tol}$     \\
			& 0 equalities             & Presented solver ($F_2$)          &  13,969.50                & 49,267,409    &  $< \epsilon_\text{tol}$     \\
			& 0 inequalities            & Presented solver ($F_3$)     &  0.91         &  789  &    $< \epsilon_\text{tol}$    \\
			&              & Presented solver ($F_4$) &   11.50               &  11,021     &  $< \epsilon_\text{tol}$   \\
			&              & Presented solver (envelope) &   0.11               &  57     & $< \epsilon_\text{tol}$   \\
			\hline
			\multicolumn{6}{l}{variable overflow: Crashed due to variable overflow in FILIB\texttt{++}.} \\
		\end{tabular}
	}
\end{table}
\subsection{Compressor Plant}
\label{subsec:Numerical Example 3}
\begin{figure}[ht]
	\begin{center}
		\includegraphics[width=4.0 in]{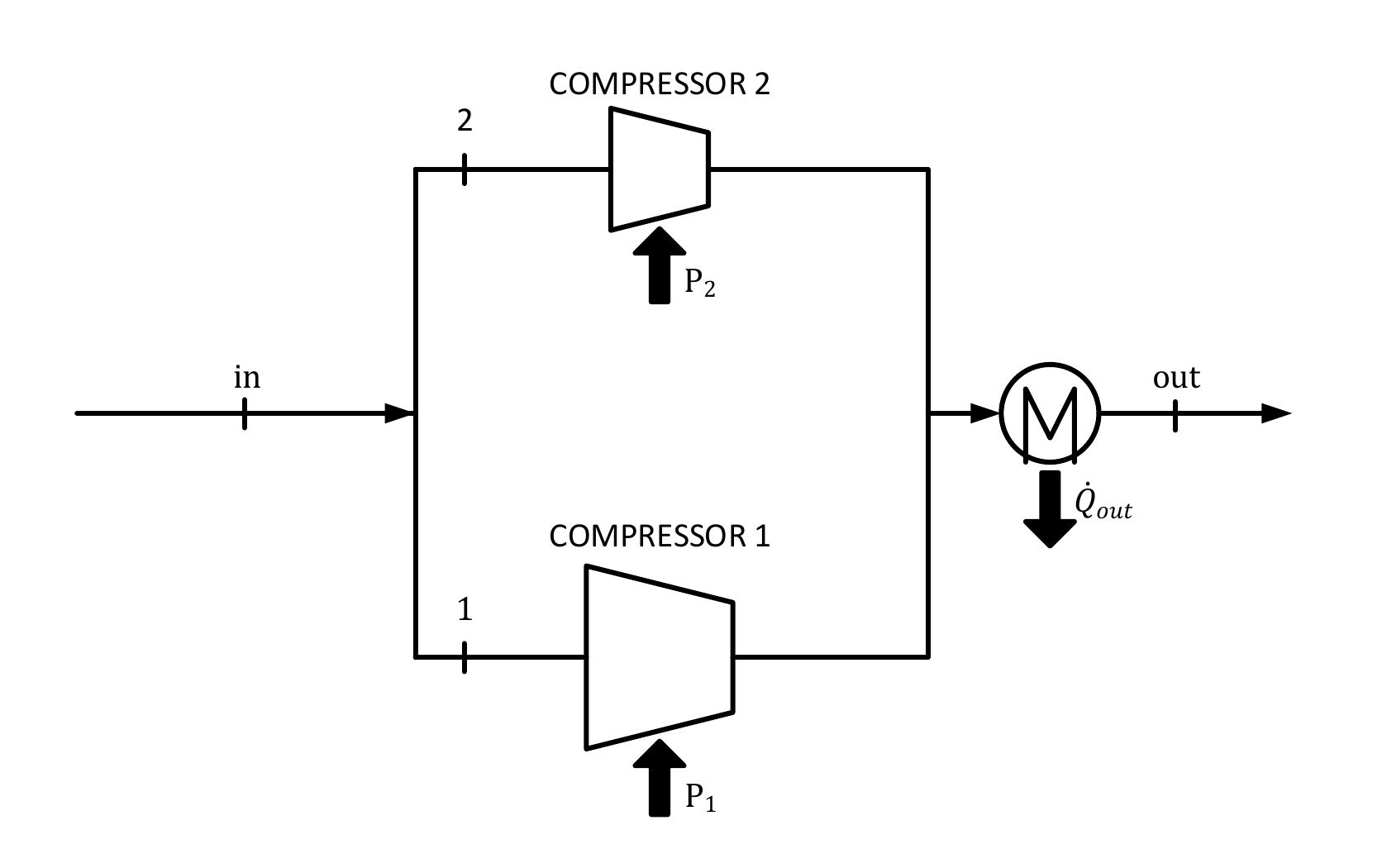}
	\end{center}
	\caption{Graphical illustration of the compressor configuration}
	\label{fig:4} 
\end{figure}

\begin{figure}[ht]
	\begin{subfigure}[c]{0.5\textwidth}
		\includegraphics[width=2.0 in]{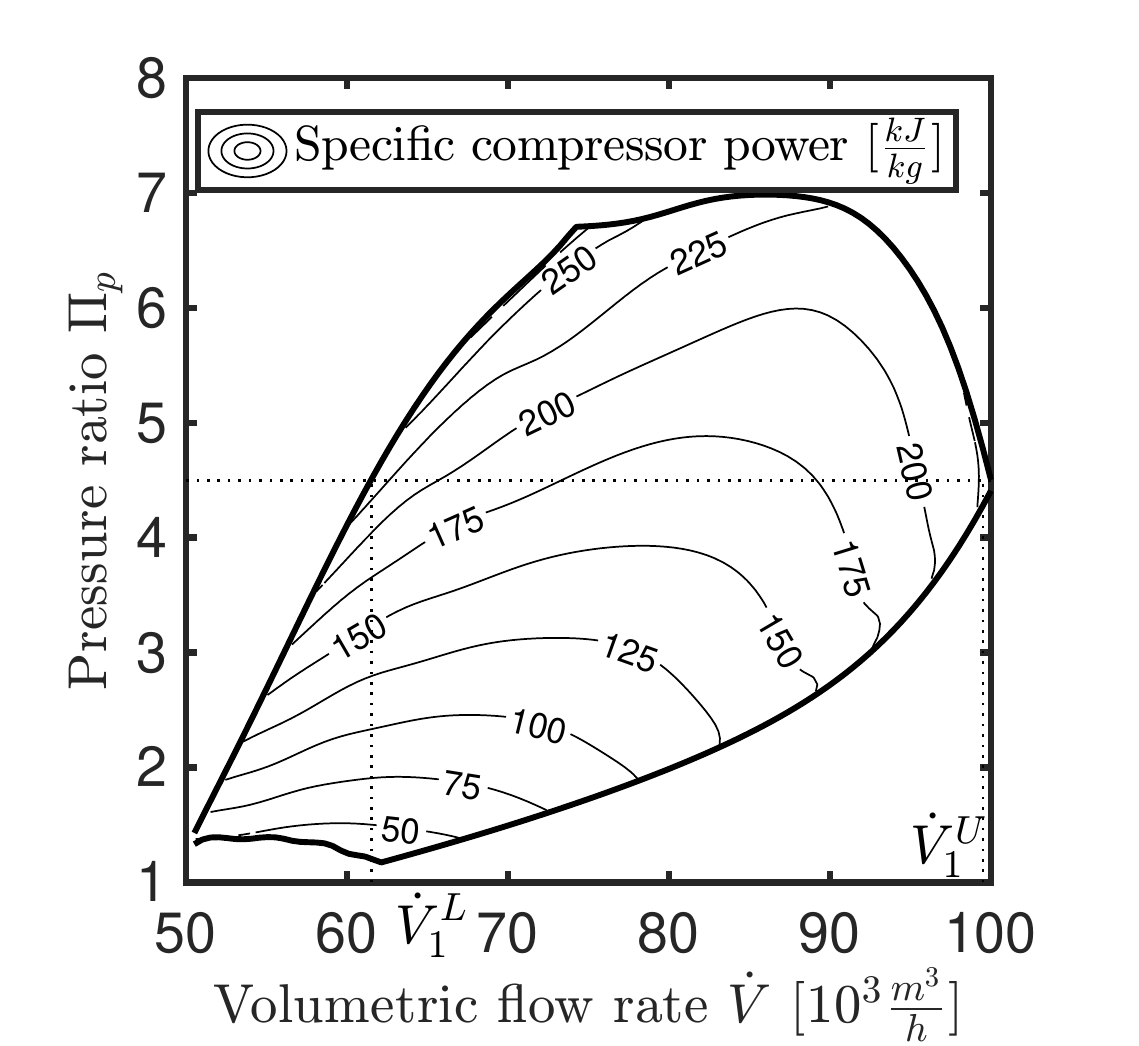}
		\subcaption{Performance map of compressor 1}
	\end{subfigure}
	\begin{subfigure}[c]{0.5\textwidth}
		\includegraphics[width=2.0 in]{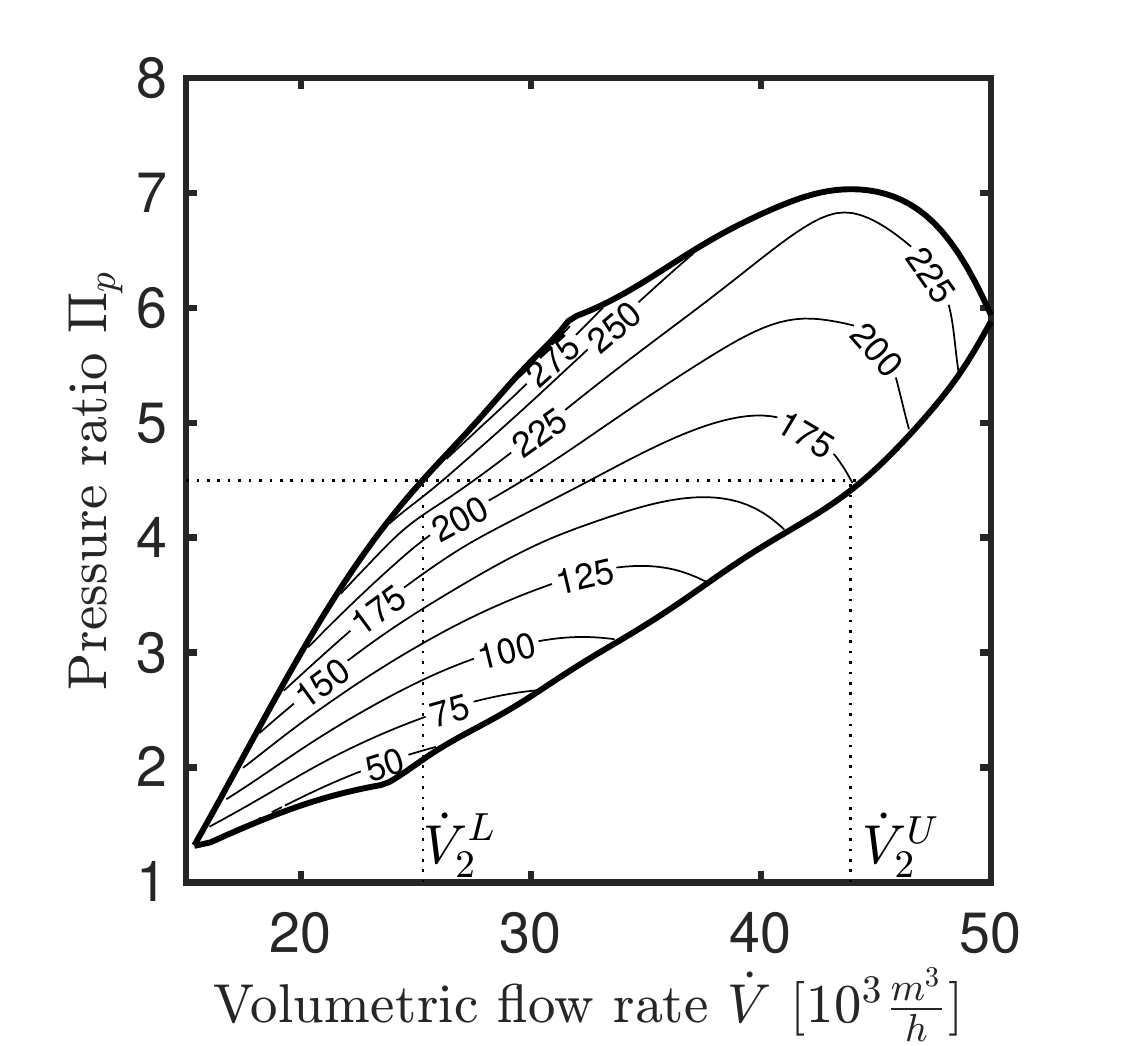}
		\subcaption{Performance map of compressor 2}
	\end{subfigure}
	\caption{Graphical illustration of the compressors' specific performance maps}
	\label{fig:3}
\end{figure}

In the third numerical example, the operating point of a compressor plant is optimized. This is a relevant case study because air compressors are commonly used in industry and have a high electrical power consumption, e.g., in cryogenic air separation units. In addition, the power consumption of industrial compressors is often provided in form of data (compressor maps, e.g., Figure \ref{fig:3}) and not mechanistic models preventing model-based optimization techniques. The considered compressor plant is comprised of two compressors that are connected in parallel as shown in Figure \ref{fig:4}. The compressors are sized such that a large compressor (compressor 1) is supplemented by a smaller compressor (compressor 2). The intent is to minimize the electrical power consumption of the overall process by optimal operation. Mathematically, the case study is challenging as it combines two MLPs with two hidden layers in one optimization problem. 
 \newline \indent
The specific power of the compressors is given by compressor maps $M_i:\mathbb{R}^2 \to \mathbb{R}$ (Figure \ref{fig:3}) that map the volumetric inlet flow rate, $\dot{V_i}$, and compression ratio, $\Pi_i = \frac{p_{\text{out},i}}{p_{\text{in},i}}$, to the specific power, $w_i$, of the corresponding compressor $i$ (maps were altered from industrial compressor maps). For the MLP training, a set of 405 data points is read out from each compressor map and randomly divided into a training, validation and test set with the respective size ratios of $0.7:0.15:0.15$. Ghorbanian and Gholamrezaei compared different regression models for compressor performance prediction and concluded that the MLP is the most powerful candidate \cite{Ghorbanian.2009}. Thus, this numerical example utilizes the network structure suggested by Ghorbanian and Gholamrezaei, namely a MLP with two hidden layers with 10 neurons each. For training, Levenberg-Marquardt algorithm and early stopping were used on a scaled training data set. The training results in a MSE of 0.17 (0.95), 0.78 (3.87) and 0.39 (1.67) on the training, validation and test set respectively for compressor one (two). 
 \newline \indent
The optimization problem minimizes the total power consumption of the compressor plant, $P_{\text{total}}$. It has one degree of freedom, the split factor ($x$), that determines the ratio between the volumetric flow rate $\dot{V}_1$ and the total volumetric inlet flow rate, $\dot{V}_{\text{in}}$. The problem is formulated as follows:
\begin{equation}
\label{eqn:CaseStudy3_OptProblem}
\begin{aligned}
&\underset{0 \leq x \leq 1}{\text{min}}  &	P_{\text{total}}(x) = &	w_{1,\text{MLP}}(\Pi_p, \dot{V}_{1}  ) \cdot \frac{{\dot{V}_{\text{in}} \cdot x}}{v_{\text{in}}} + w_{2,\text{MLP}}(\Pi_p, \dot{V}_{2}  ) \cdot \frac{\dot{V}_{\text{in}} \cdot (1-x)}{v_{\text{in}}} \\
& \text{s.t.} & \dot{V}_{1} - \dot{V}^{U}_{1} & \le 0 \\
& 			  & \dot{V}^{L}_{1} - \dot{V}_{1} & \le 0 \\
&  			  & \dot{V}_{2} - \dot{V}^{U}_{2} & \le 0 \\
& 			  & \dot{V}^{L}_{2} - \dot{V}_{2} & \le 0 \\
\end{aligned}
\end{equation}
with the constants $v_{\text{in}} = 0.8305 ~ \frac{\text{m}^3}{\text{kg}}, ~ p_{\text{in}} = 100 ~ \text{kPa}, ~ p_{\text{out}} = 450 ~ \text{kPa}, ~ \dot{V}_{\text{in}} = 100 \cdot 10^{3} ~  \frac{\text{m}^3}{\text{h}}, ~ \dot{V}^{L}_{1} = 61.5 \cdot 10^{3} ~ \frac{\text{m}^3}{\text{h}}, ~ \dot{V}^{U}_{1} = 100 \cdot 10^{3} ~ \frac{\text{m}^3}{\text{h}}, ~ \dot{V}^{L}_{2}= 25.3 \cdot 10^{3} ~ \frac{\text{m}^3}{\text{h}} ~ \text{and} ~ \dot{V}^{U}_{2} = 44.4 \cdot 10^{3} ~  \frac{\text{m}^3}{\text{h}}$. Further, $\dot{V}_{1} = \dot{V}_{in} \cdot x ~ \text{and} ~ \dot{V}_{2} = \dot{V}_{\text{in}} \cdot (1-x)$ are explicit functions of $x$. The functions $w_{1,\text{MLP}}$ and $w_{2,\text{MLP}}$ describe the specific power predicted by the MLPs. The lower and upper bounds of the volumetric flow rates, $\dot{V}^{L}_{1}$, $\dot{V}^{L}_{2}$ and $\dot{V}^{U}_{1}$, $\dot{V}^{U}_{2}$ correspond to the surge and choke lines of the compressors at the specified pressure ratio of 4.5 (left and right hand boundaries in Figure \ref{fig:3}). \newline \indent 
The optimization problem \eqref{eqn:CaseStudy3_OptProblem} is solved using an absolute optimality gap of $10^{-4}$ and a relative optimality gap of $10^{-12}$. Further, a CPU time limit of 100,000 seconds is set. As shown in Table \ref{tab:NumericalResults_Compressor}, the RS formulation reduces the number of variables from 97 to 1 and eliminates all 96 equality constrains compared to the FS formulation. The RS formulation is solved using the presented solver within 0.23 seconds when using the envelopes. This is over four hundred thousand times faster than the solution with the same solver and relaxation in the FS formulation. BARON converges only using reformulation $F_2$ (within 25.02 second). When the other reformulations are used, the desired tolerance is not reached within 100,000 seconds. Thus, the proposed approach is over one hundred times faster than BARON if BARON uses $F_2$ and over four hundred thousand times faster if BARON uses $F_1$, $F_3$ or $F_4$.
\begin{table}[ht]
	\centering
	\caption{Numerical results of the compressor plant optimization (Subsection \ref{subsec:Numerical Example 3})}
	\label{tab:NumericalResults_Compressor}
	\resizebox{\columnwidth}{!}{
		\begin{tabular}{lllrrr}
			\hline
			& Problem size & Solver                      & CPU time {[}s{]} & Iterations 				& Abs. gap \\ \hline
			\multirow{4}{*}{(FS)} 
			& 97 variables  & BARON ($F_1$)                &   100,000.00       &  1,934,517 		    &  $> 3 \cdot 10^8$     \\
			& 96 equalities & BARON ($F_2$)                &   25.02         &  892        		&  $< \epsilon_\text{tol}$\\
			& 4 inequalities             & BARON ($F_3$)   &  100,000.00     &  4,224,962    		&     0.6\\
			&              & BARON ($F_4$)                &   100,000.00      & 3,128,979 &    0.8 \\
			
			&  				& Presented solver ($F_1$)      &  	\multicolumn{3}{c}{--variable overflow--}   \\
			&             	& Presented solver ($F_2$)      &  	\multicolumn{3}{c}{--variable overflow--}    \\
			&             	& Presented solver ($F_3$)     	& 	\multicolumn{3}{c}{--variable overflow--}     \\
			&              	& Presented solver ($F_4$) 		&  	\multicolumn{3}{c}{--variable overflow--}   \\
			&              	& Presented solver (envelope) 	&    100,000.00            &   846,124     & 181.2   \\
			\hline
			\multirow{5}{*}{(RS)}  
			& 1 variables & Presented solver ($F_1$)          &  \multicolumn{3}{c}{--variable overflow--} \\
			& 0 equalities             & Presented solver ($F_2$)          &  \multicolumn{3}{c}{--variable overflow--}    \\
			& 4 inequalities            & Presented solver ($F_3$)     &  0.34         &   171  &   $< \epsilon_\text{tol}$    \\
			&              & Presented solver ($F_4$) &    \multicolumn{3}{c}{--variable overflow--} \\
			&              & Presented solver (envelope) &   0.23             &  81    &   $< \epsilon_\text{tol}$  \\
			\hline
			\multicolumn{6}{l}{variable overflow: Crashed due to variable overflow in FILIB\texttt{++}.} \\
		\end{tabular}
	}
\end{table}  \newline \indent
The optimal solution of the problem is $P_{\text{total}}=6.20 ~ \text{MW}$ with a split ratio of $x = 0.683$ that corresponds to $\dot{V}_1 = 68.28 ~ \frac{m^3}{h}$ and $\dot{V}_2 = 31.72 ~ \frac{m^3}{h}$. The optimization using learned compressor maps can have some advantage over simple operating heuristics. For instance, if the main air compressor would be operated at its most efficient operating point, the the total power consumption would be 6.2\% higher ($P_{\text{total}}=6.61 ~ \text{MW} $) and if the flowrates would be divided according to the maximum volumetric capacity of the compressors, the total power consumption would be 9.3\% higher ($P_{\text{total}}=6.84 ~ \text{MW}$).
\subsection{Cumene Process}
\label{subsec:Numerical Example 4}
\begin{figure}[ht]
	\begin{center}
		\includegraphics[width=4.0 in]{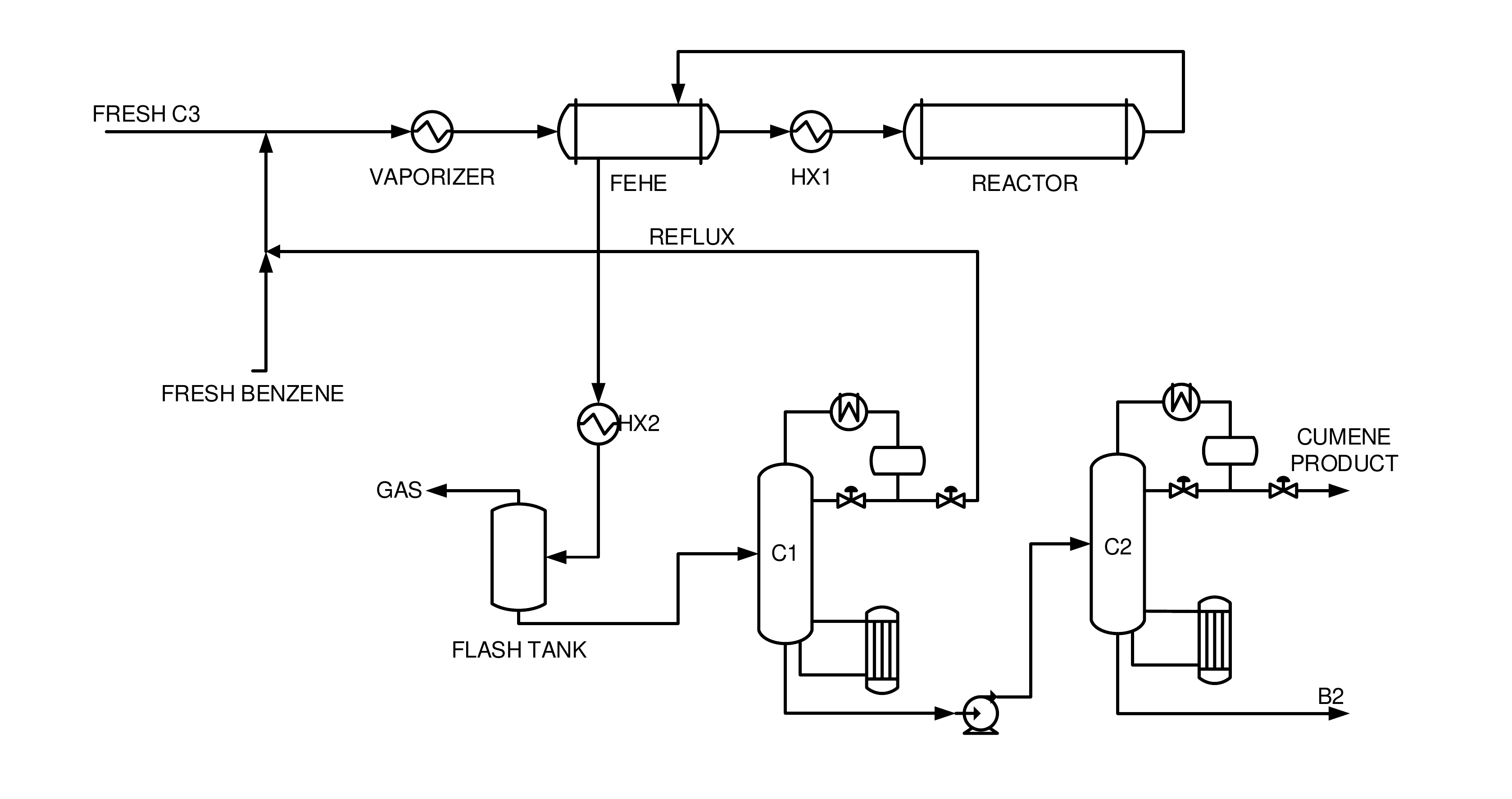}
	\end{center}
	\caption{Graphical illustration of the cumene process}
	\label{fig:5}
\end{figure}

\begin{figure}[ht]
	\begin{center}
		\includegraphics[width=4.0 in]{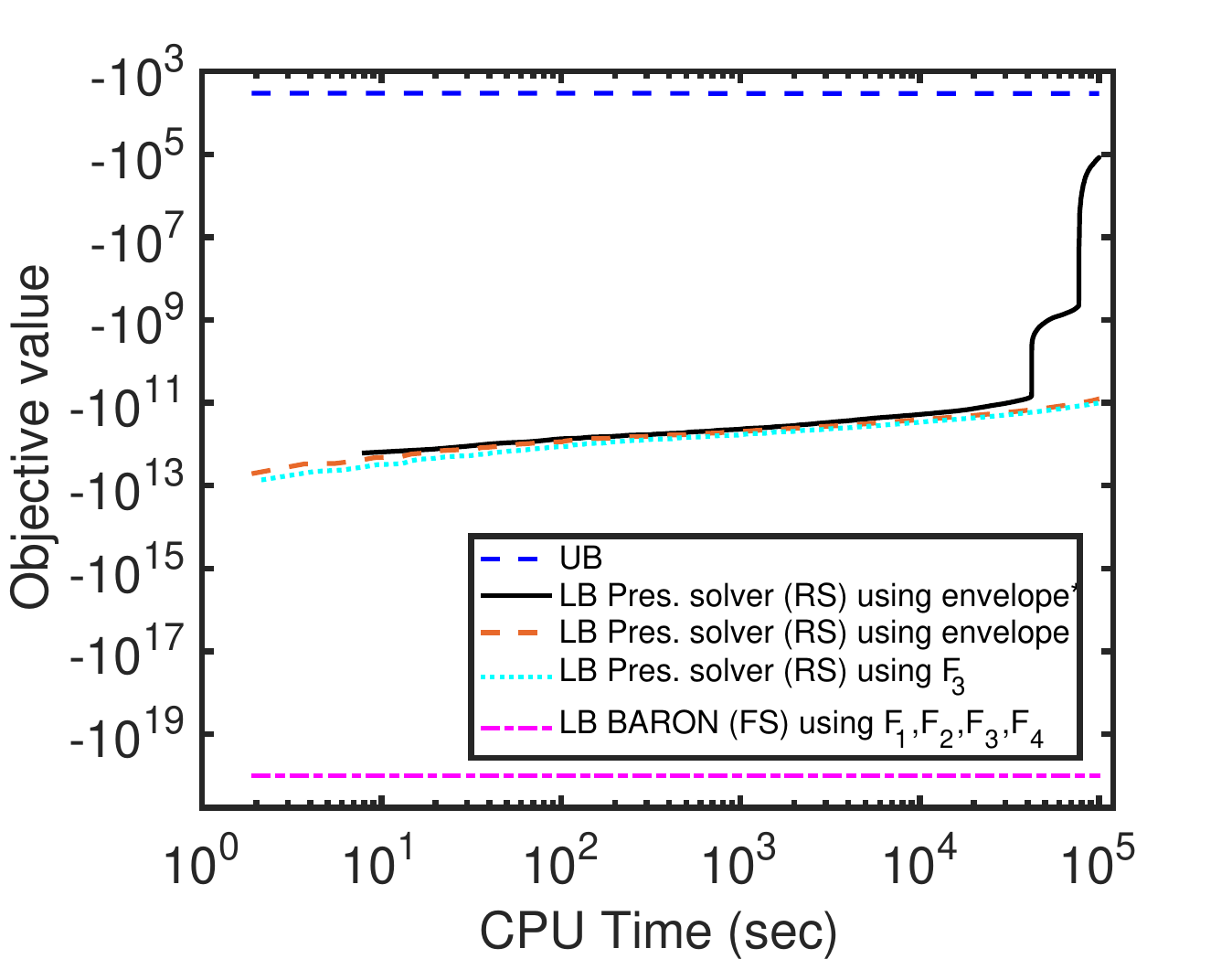}
	\end{center}
	\caption{Convergence plot of the cumene process optimization problem (Subsection \ref{subsec:Numerical Example 4}). Herein, UB refers to the upper bound and LB to the lower bound of the corresponding solvers. The asterisk (*) refers to adapted solver options such that the upper bound is obtained by function evaluation, not local search}
	\label{fig:6} 
\end{figure}

In the third numerical example, the operating point of a cumene process is optimized illustrating a complex industrial unit. As illustrated in Figure \ref{fig:5}, the cumene process consists of a plug flow reactor, two rectification columns, one flash, several (integrated) heat exchangers and recycles. A detailed process description including necessary details for model building can be found in \cite{Luyben.2010}. For this work, we optimize a hybrid model consisting of 14 MLPs that emulate the process. This hybrid model was provided by Schultz \textit{et al.} \cite{Schultz.2016} who modeled the cumene process in ASPEN Plus and learned the MLPs using simulated data. The optimization problem minimizes the negative total profit $P$ of the process \cite{Schultz.2016}:
\begin{equation}
\label{eqn:CaseStudy4_OptProblem}
\begin{aligned}
&\underset{\textbf{x}}{\text{min}}  & & -P =-  F_{\text{Prod}} \cdot \$_{\text{Cumpure}} + Fresh_{\text{C3}} \cdot \$_{\text{C3}} + Fresh_{\text{benz}} \cdot \$_{\text{benzpure}} \\
& & & - F_{\text{Gas}} \cdot (  Zgas_{\text{benz}} \cdot \$_{\text{benz}} + Zgas_{\text{cum}} \cdot \$_{\text{cum}} \\
& & & + Zgas_{\text{propy}} \cdot \$_{\text{propy}} + Zgas_{\text{propa}} \cdot \$_{\text{propa}} \\
& & & + Zgas_{\text{PBI}} \cdot \$_{\text{PBI}} )  - F_{\text{B2}} \cdot ( ZB2_{\text{cum}} \cdot \$_{\text{cum}} + ZB2_{\text{PBI}} \cdot \$_{\text{PBI}}) \\
& & & - Q_{\text{Reactor}} \cdot \$_{\text{Eger}} + (Q_{\text{boiler}} + Q_{\text{rebC1}} + Q_{\text{rebC2}} )\cdot \$_{\text{HP}}\\
& & & + (Q_{\text{HX1}} + Q_{\text{HX2}})\cdot \$_{\text{EE}}\\
& \text{s.t.}& & T_{\text{reactor}} - 390^\circ \text{C}   \le 0, \quad 360^\circ \text{°C} - T_{\text{reactor}} \le 0,\\
&  			 & & Q_{\text{rebC1}} - 1.97 \text{Gcal/h}   \le 0, \quad 1.31 \text{Gcal/h}- Q_{\text{rebC1}} \le 0,\\
&  			 & & RR_{\text{C1}} - 1   \le 0, \quad 0.37 - RR_{\text{C1}} \le 0,\\
&  			 & & DF_{\text{C2}} - 0.99  \le 0, \quad 0.85 - DF_{\text{C2}} \le 0,\\
&  			 & & RR_{\text{C2}} - 1.2   \le 0, \quad 0.2 - RR_{\text{C2}} \le 0,\\
&			 & &  0.999 - Zprod_{\text{cum}}  \le 0
\end{aligned}
\end{equation}
with $\textbf{x} = \left( T_{\text{reactor}}, Q_{\text{rebC1}}, RR_{\text{C1}}, DF_{\text{C2}},RR_{\text{C2}} \right)$. A list of all economic parameters ($\$_i$) of the optimization problem can be found in \cite{Schultz.2016}.
 \newline \indent
The hybrid model is composed by 14 MLPs, each consisting 10 to 20 neurons, resulting in a FS optimization problem with 794 variables, 789 equality and 1 inequality constraints. The RS formulation constitutes only 5 variables, 0 equality and 1 inequality constraints. Due to the complexity of the case study, the problem is only implemented in the FS formulation in BARON and the RS formulation in the presented solver. The FS formulation in the presented solver is omitted as it is always outperformed by the RS formulation in the previous examples. 
\begin{table}[ht]
	\centering
	\caption{Numerical results of the cumene process optimization problem (Subsection \ref{subsec:Numerical Example 4})}
	\label{tab:NumericalResults_Cumene}
	\resizebox{\columnwidth}{!}{
		\begin{tabular}{lllrrr}
			\hline
			& Problem size & Solver                      & CPU time {[}s{]} & Iterations 				& Abs. gap \\ \hline
			\multirow{4}{*}{(FS)} 
			& 794 variables  & BARON ($F_1$)                & 100,000.00       	&  110,453    	&  $ 1 \cdot 10^{20}$    \\
			& 789 equalities & BARON ($F_2$)                & 100,000.00        &  120,536 		&  $ 1 \cdot 10^{20}$   \\
			& 1 inequalities& BARON ($F_3$)   			  	& 100,000.00     	&  136,645		&  $ 1 \cdot 10^{20}$     \\
			&              & BARON ($F_4$)                	& 100,000.00        &  117,400		&  $ 1 \cdot 10^{20}$   \\
			\hline
			\multirow{5}{*}{(RS)}  
			& 5 variables & Presented solver ($F_1$)        &  \multicolumn{3}{c}{--variable overflow--}   \\
			& 0 equalities & Presented solver ($F_2$)       &  \multicolumn{3}{c}{--variable overflow--}     \\
			& 1 inequalities & Presented solver ($F_3$)     &   100,000.00        &   4,772,133   & $1 \cdot 10^{11}$   \\
			&              & Presented solver ($F_4$)		&  \multicolumn{3}{c}{--variable overflow--}   \\
			&              & Presented solver (envelope) 	&   100,000.00        &   5,683,103   & $8 \cdot 10^{10}$   \\
			&              & Presented solver (envelope)* 	&   100,000.00        &   12,939,508  & $1 \cdot 10^5 ~$   \\
			\hline
			\multicolumn{6}{l}{variable overflow: Crashed due to variable overflow in FILIB\texttt{++}.} \\
			\multicolumn{6}{l}{*: Adapted setting such that upper bound is obtained by function evaluation, not local search.} \\
		\end{tabular}
	}
\end{table} \newline \indent
The results in Table \ref{tab:NumericalResults_Cumene} show that none of the tested solution approaches converge to the desired tolerance within the CPU time limit. In order to analyze this, a  convergence plot is provided in Figure \ref{fig:6} that depicts the lower bounds of the solvers over CPU time. Apparently, BARON does not improve its initial lower bound on the objective at all. In comparison, the presented solver improves its lower bound steadily but slowly. Due to the high complexity of the case-study another optimization run (annotated with the asterisk (*)) is executed using adapted solver options. More precisely, instead of using a local NLP solver in each iteration for upper bounding, the objective function is just evaluated at the center of the current interval. This leads to about 2.3 times more iterations in the same CPU time and further improvements of the lower bound. However, the complex example would still require longer CPU times to converge to the desired tolerance.
\section{Conclusion and Future Work}
\label{sec:Conclusionandfuturework}
A method for deterministic global optimization of ANN embedded optimization problems is presented. The proposed method formulates the problem in a RS and computes lower bounds by propagation of McCormick relaxations through the network equations. Thus, variables inside the ANNs and the network equations are not visible to the optimization algorithm and problem size is reduced significantly. \newline \indent 
The proposed method is tested on four numerical examples including one illustrative function and three engineering applications.
In all examples, the number of optimization variables and equality constraints could reduced by the RS formulations. In particular, problems including ANNs with a large number of neurons benefit from the RS formulation. Further, the RS formulation results in all examples in an acceleration of the optimization by factors of up to four hundred thousand compared to the FS formulation using the same solver. A direct comparison to the state-of-the-art solver BARON, which works in FS, is more difficult because BARON currently does not include the hyperbolic tangent activations function. However, the results suggest that the combination of the RS formulation and the envelopes yield favorable performance of the proposed method in comparison to BARON. \newline \indent
In terms of problem sizes, the proposed approach optimized single shallow ANNs with up to 700 neurons in under 10 seconds. An analysis of the scaling of the approach also reveals that deep network require more CPU time for optimization than shallow networks with the same or even larger number of neurons. Thus, we conclude that shallow ANNs are currently more suitable for efficient optimization than deep ANNs. The engineering examples further embed up to 14 ANNs. Although, two of the engineering applications could be solved within a fraction of a second, the solution of the complex cumene process case-study converges slowly. This is probably due to the propagation of relaxations through a large number of ANNs and equations. In general, the considered problem sizes significantly exceed the examples for deterministic ANN embedded optimization found in the literature (i.e., Smith \textit{et al.} (2013) \cite{Smith.2013} optimized an ANN with 1 layer, 3 neurons and 5 inputs using BARON).\newline \indent
As shown in the case studies, the combination of ANNs as universal approximators and an efficient deterministic global optimization algorithm is a powerful tool for various applications. In addition, this method can be further extended to, for instance, mixed integer problems where a process superstructure is optimized globally \cite{Henao.2011,Lee.2016}. Here, unit operations, thermodynamic models and even dynamics could be emulated by ANNs allowing the utilization of data from different sources such as experiments, process simulations and life-cycle assessment tools (e.g., \cite{Helmdach.2017}). Also, complex model parts could be replaced by ANNs that lump these complex parts yielding possibly tight relaxations and less optimization variables. Other possible extension of the proposed method could be more efficient methods for optimization of deep ANNs with a large number of neurons as well as deterministic global training of ANNs.

\section{Acknowledgements}
The authors gratefully acknowledge the financial support of the Kopernikus-project SynErgie by the Federal Ministry of Education and Research (BMBF) and the project supervision by the project management organization Projekttr\"ager J\"ulich (PtJ). We are grateful to Jaromi{\l} Najam, Dominik Bongartz and Susanne Scholl for their work on MAiNGO and Beno{\^{i}}t Chachuat for providing MC\texttt{++}. We thank Eduardo Schutz for providing the model of the fourth numerical example, Adrian Caspari \& Pascal Sch\"afer for helpful discussions and Linus Netze \& Nils Gra\ss ~ for implementing case studies.

\appendix

\section{Appendix}
\subsection{Convex and Concave Envelopes of the Hyperbolic Tangent Activation Function}
\label{sec:Convex and concave envelopes of the hyperbolic tangent activation function}
In this subsection, the envelopes of the hyperbolic tangent ($\tanh$) function are derived on a compact interval $D = [x^{L}, x^{U}]$. As the hyperbolic tangent function is one-dimensional, McCormick \cite{McCormick.1976}  gives a method to construct its envelopes. More specifically, as the hyperbolic tangent function is convex on $(-\infty,0]$ and concave on $[0,+\infty)$, its convex envelope, $F^{cv}: \mathbb{R} \to \mathbb{R}$, and concave envelope, $F^{cc}: \mathbb{R} \to \mathbb{R}$, are given:
\begin{equation}\label{eqn:funderlin1}
F^{cv}(x)= 
\begin{cases}
\tanh(x), &\quad  x^{U} \le 0 \\
\sec(x),  &\quad 0 \le x^{L} \\
F^{cv}_3(x),  &\quad \text{otherwise}
\end{cases} 
\end{equation}
\begin{equation}\label{eqn:funderlin2}
F^{cc}(x)= 
\begin{cases}
\sec(x),  &\quad x^{U} \le 0   \\
\tanh(x),  &\quad 0 \le x^{L} \\
F^{cc}_3(x),  &\quad \text{otherwise}
\end{cases} 
\end{equation}
where the secant given as $\sec(x)=\frac{\tanh(x^{U}) - \tanh(x^{L})}{x^{U} - x^{L}} x + \frac{x^{U} \tanh(x^{L}) - x^{L} \tanh(x^{U})}{x^{U} - x^{L}}$. For $x^L < 0 < x^U $, the hyperbolic tangent function is nonconvex and nonconcave. The convex envelope, $F^{cv}_3: \mathbb{R} \to \mathbb{R}$, for this case is:
\begin{equation}\label{eqn:Fcv_3}
F^{cv}_3(x)= 
\begin{cases}
\tanh(x), &\quad x \le x^u_c \\
\frac{\tanh(x^{U})-\tanh(x^u_c)}{x^{U}-x^u_c} \cdot (x-x^u_c) + \tanh(x^u_c),  &\quad x > x^u_c
\end{cases} 
\end{equation}	
where $x^u_c = \max(x^{u*}_c,x^{L})$ and $x^{u*}_c$ is the solution of:
\begin{equation}\label{eqn:x_c^u}
1 - \tanh^2(x^u_c)=\frac{\tanh(x^{U})-\tanh(x^u_c)}{x^{U}-x^u_c} , \quad x^u_c \le 0
\end{equation}
solved numerically for every interval. Similarly, the concave envelope, $F^{cc}_3: \mathbb{R} \to \mathbb{R}$, of $\tanh$ is obtained for $x^L < 0 < x^U$ as:
\begin{equation}\label{eqn:Fcc_3}
F^{cc}_3(x)= 
\begin{cases}		
\frac{\tanh(x^o_c)- \tanh(x^{L})}{x^o_c-x^{L}} \cdot (x-x^{L}) + \tanh(x^{L}),  &\quad x < x^o_c \\
\tanh(x), &\quad x \ge x^o_c
\end{cases} 
\end{equation}	
where $x^o_c = \min(x^{o*}_c,x^{U})$ and $x^{o*}_c$ is the solution of:
\begin{equation}\label{eqn:x_c^o}
1 - \tanh^2(x^o_c)=\frac{\tanh(x^o_c)-\tanh(x^{L})}{x^o_c-x^{L}} , \quad x^o_c \ge 0
\end{equation}	\indent
In the following, we show that the convex and concave envelopes of the hyperbolic tangent function are smooth ($C^1$) and strictly monotonically increasing.
\begin{proposition} \label{prop:smoothness of hyperbolic tangent relaxations}
	(smoothness of hyperbolic tangent relaxations).  The convex and concave envelopes of the hyperbolic tangent function, $F^{cv}(x)$ and $F^{cc}$, are once continuously differentiable ($C^1$) and in general not $C^2$.
\end{proposition}
{\it Proof}
	For $x^U \le 0$, $F^{cv}(x) = \tanh(x)$ and $F^{cc}(x)=\sec(x)$ which are $C^\infty$. Similarly, for $0 \le x^L$, $F^{cc}(x) = \tanh(x)$ and $F^{cv}(x)=\sec(x)$ which are $C^\infty$. For  $x^L < 0 < x^U$, the envelopes are given by \eqref{eqn:Fcv_3} and \eqref{eqn:Fcc_3}. These are at least once continuously differentiable ($C^1$) because of \eqref{eqn:x_c^u} and \eqref{eqn:x_c^o}. \eqref{eqn:Fcv_3} and \eqref{eqn:Fcc_3} are at most $C^1$ because 	
	\begin{equation}\label{eqn:d2/dx2_F3_cv}
	\left. \frac{ \text{d}^2(F_3^{cv})}{\text{d} x^2} \right|_{x} = 
	\begin{cases}		
	-2 \sech^2(x) \tanh(x),  &\quad x < x^u_c \\
	0, &\quad x \ge x^u_c
	\end{cases} 
	\end{equation}	
	and
	\begin{equation}\label{eqn:d2/dx2_F3_cc}
	\left. \frac{ \text{d}^2(F_3^{cc})}{\text{d} x^2} \right|_{x} = 
	\begin{cases}		
	0,  &\quad x < x^o_c \\
	-2 \sech^2(x) \tanh(x), &\quad x \ge x^o_c
	\end{cases} 
	\end{equation}	
	where $\sech(x)$ is the hyperbolic secant function,  are not continuous at $x^u_c$ and $x^o_c$, respectively.
\qed
As shown in Proof \ref{prop:smoothness of hyperbolic tangent relaxations}, the first derivative of the convex and concave envelopes of the hyperbolic tangent function, $\frac{ \text{d}(F^{cv})}{\text{d} x}$ and $\frac{ \text{d}(F^{cc})}{\text{d} x}$, are continuous but not continuously differentiable. The following proof shows that the first derivative of the convex and concave envelopes of the hyperbolic tangent function are Lipschitz continuous.
\begin{proposition} \label{prop:lipschitz}
	(Lipschitz continuity of 1st derivative of hyperbolic tangent relaxations). The
	second derivative of the convex and concave envelopes of the hyperbolic tangent function are bounded. This implies that the first derivative of the convex and concave envelopes of the hyperbolic tangent function, $\frac{ \text{d}(F_3^{cv})}{\text{d} x}$ and $\frac{ \text{d}(F_3^{cc})}{\text{d} x}$, are at least once Lipschitz continuous.
\end{proposition}
{\it Proof}
	For $x^U \le 0$, $F^{cv}(x) = \tanh(x)$ and $F^{cc}(x)=\sec(x)$ which are $C^\infty$. Similarly, for $0 \le x^L$, $F^{cc}(x) = \tanh(x)$ and $F^{cv}(x)=\sec(x)$ which are $C^\infty$. For  $x^L < 0 < x^U$, the second derivative of the convex and concave envelopes of the hyperbolic tangent function are given by \eqref{eqn:d2/dx2_F3_cv} and \eqref{eqn:d2/dx2_F3_cc}. This implies that
	\begin{equation}\label{eqn:proofL1}
	\left| \frac{ \text{d}^2(F_3^{cv}(x))}{\text{d} x^2} \right| \le 2 \left|\sech^2(\tilde{x}^{cv}) \tanh(\tilde{x}^{cv}) \right| =  2\left| \frac{\sinh(\tilde{x}^{cv})}{\cosh^3(\tilde{x}^{cv})} \right|
	\end{equation}
	with $x^L \le x \le x^U$ and $x^L \le \tilde{x}^{cv} \le x_c^u$. From $\cosh(x) \ge 1$ it follows that 
	\begin{equation}\label{eqn:proofL2}
	2\left| \frac{\sinh(\tilde{x}^{cv})}{\cosh^3(\tilde{x}^{cv})} \right| \le 2\left| \sinh(\tilde{x}^{cv}) \right|
	\end{equation}
	As $\sinh(x)$ is a monotonic function that is point symmetric with respect to the origin and $x^L \le \tilde{x}^{cv} \le 0$, it follows that the second derivative of the convey envelope is bounded
	\begin{equation}\label{eqn:proofL3}
	2\left| \sinh(\tilde{x}^{cv}) \right| \le 2 \sinh(\left| x^L \right|)
	\end{equation}
	Thus, the first derivative of the convex envelope of the hyperbolic tangent function, $\frac{ \text{d}(F^{cv})}{\text{d} x}$, is at least once Lipschitz continuous with a Lipschitz constant of at most $L^{cv} = 2 \sinh(\left| x^L \right|)$. \newline			
	Similarly, it holds that
	\begin{equation}\label{eqn:proofL4}
	\left| \frac{ \text{d}^2(F_3^{cc}(x))}{\text{d} x^2} \right| \le 2 \left|\sech^2(\tilde{x}^{cc}) \tanh(\tilde{x}^{cc}) \right| \le  2\left| \sinh(\tilde{x}^{cc}) \right|
	\end{equation}
	with $x^L \le x \le x^U$ and $0 \le x_c^o \le \tilde{x}^{cc} \le x^U$. It follows that, the first derivative of the concave envelope of the hyperbolic tangent function, $\frac{ \text{d}(F^{cc})}{\text{d} x}$, is at least once Lipschitz continuous with a Lipschitz constant of at most $L^{cc} = 2 \sinh(\left| x^U \right|)$. 
\qed
\begin{proposition} \label{prop:monotonicity of hyperbolic tangent relaxations}
	(monotonicity of hyperbolic tangent relaxations).   The convex and concave envelopes of the hyperbolic tangent function are strictly monotonically increasing.
\end{proposition}
{\it Proof}
	For $x^U \le 0$, $F^{cv}(x) = \tanh(x)$ and $F^{cc}(x)=\sec(x)$. Similarly, for $0 \le x^L$, $F^{cc}(x) = \tanh(x)$ and $F^{cv}(x)=\sec(x)$. As $\frac{ \text{d}(\tanh(x))}{\text{d} x} = 1 - \tanh^2(x) > 0$, $\tanh(x)$ is strictly monotonic monotonically. As $x^U > x^L$, $\sec(x)$ is strictly monotonically increasing. For  $x^L < 0 < x^U$, the envelopes are given by \eqref{eqn:Fcv_3} and \eqref{eqn:Fcc_3}. These are again strictly monotonically increasing because $x_c^o > x^L$ and $x^U > x_c^u$.
\qed
\subsection{Convex and Concave Envelopes of the Sigmoid Activation Function}
\label{sec:Convex and concave envelopes of the sigmoid activation function}
Another common activation function of ANNs is the sigmoid function. The sigmoid function can be reformulated using $\text{sig}(x) = \frac{1}{2}(1+\tanh(\frac{x}{2}))$. The convex and concave envelopes of the sigmoid function, $F^{cv}_{\text{sig}}(x)$ and $F^{cc}_{\text{sig}}(x)$, can be derived using the reformulation and the convex and concave envelopes of the hyperbolic tangent function, $F^{cv}_{\tanh}(x)$ and $F^{cc}_{\tanh}(x)$:
\begin{equation}\label{eqn:cv_sig}
F^{cv}_{\text{sig}}(x) = \frac{1}{2}\left(1+F^{cv}_{\tanh}(\frac{x}{2})\right)
\end{equation}
\begin{equation}\label{eqn:cc_si}
F^{cc}_{\text{sig}}(x) = \frac{1}{2}\left(1+F^{cc}_{\tanh}(\frac{x}{2})\right)
\end{equation} \indent
For simplicity, this manuscript does not provide proofs for smoothness and monotonicity of the envelopes of the sigmoid function. However, similar results to the ones for the hyperbolic tangent activation function (see Appendix \ref{sec:Convex and concave envelopes of the sigmoid activation function}) can be derived for the sigmoid activation function.
\subsection{McCormick Relaxations of Reformulations of the Hyperbolic Tangent Activation Function}
\label{sec:McCormick relaxations of reformulations of the hyperbolic tangent activation function}
Reformulations of the hyperbolic tangent function are necessary for solvers where the hyperbolic tangent function is not directly available. Four common reformulations $F_1, F_2, F_3, F_4: \mathbb{R} \to \mathbb{R}$ are defined as follows: \newline
\begin{equation}
F_1(x) \coloneqq \frac{e^x - e^{-x}}{e^x + e^{-x}} =\tanh(x) 
\end{equation}

\begin{equation}
F_2(x)  \coloneqq \frac{e^{2x} - 1}{e^{2x} + 1} =\tanh(x) 
\end{equation}

\begin{equation}
F_3(x) \coloneqq 1 - \frac{2}{e^{2x} + 1} =\tanh(x) 
\end{equation}

\begin{equation}
F_4(x) \coloneqq \frac{1 - e^{-2x}}{1 + e^{-2x}}  =\tanh(x) 
\end{equation}
The convex and concave relaxations of  $F_1, F_2, F_3, F_4$ can be computed using MC\texttt{++}. 
It can be shown that the convex and concave relaxations of $F_1$, $F_2$ and $F_4$ are weaker than the ones of $F_3$ in specific cases (e.g., on the interval $x \in D=[-1,1]$). Further, the convex and concave relaxations of $F_1$, $F_2$ and $F_4$ are not differentiable. Finally, it can be observed that relaxations of the reformulations are considerably weaker than the envelopes of the hyperbolic tangent function.


\begin{thebibliography}{10}
	\providecommand{\url}[1]{{#1}}
	\providecommand{\urlprefix}{URL }
	\expandafter\ifx\csname urlstyle\endcsname\relax
	\providecommand{\doi}[1]{DOI~\discretionary{}{}{}#1}\else
	\providecommand{\doi}{DOI~\discretionary{}{}{}\begingroup
		\urlstyle{rm}\Url}\fi
	
	\bibitem{Hornik.1989}
	Hornik, K., Stinchcombe, M., White, H.: Multilayer feedforward networks are
	universal approximators.
	\newblock Neural Networks \textbf{2}(5), 359--366 (1989).
	\newblock \doi{10.1016/0893-6080(89)90020-8}
	
	\bibitem{Gasteiger.1993}
	Gasteiger, J., Zupan, J.: Neural networks in chemistry.
	\newblock Angewandte Chemie International Edition in English \textbf{32}(4),
	503--527 (1993).
	\newblock \doi{10.1002/anie.199305031}
	
	\bibitem{AzlanHussain.1999}
	{Azlan Hussain}, M.: Review of the applications of neural networks in chemical
	process control - simulation and online implementation.
	\newblock Artificial Intelligence in Engineering \textbf{13}(1), 55--68 (1999).
	\newblock \doi{10.1016/S0954-1810(98)00011-9}
	
	\bibitem{AgatonovicKustrin.2000}
	Agatonovic-Kustrin, S., Beresford, R.: Basic concepts of artificial neural
	network modeling and its application in pharmaceutical research.
	\newblock Journal of Pharmaceutical and Biomedical Analysis \textbf{22}(5),
	717--727 (2000).
	\newblock \doi{10.1016/S0731-7085(99)00272-1}
	
	\bibitem{WitekKrowiak.2014}
	Witek-Krowiak, A., Chojnacka, K., Podstawczyk, D., Dawiec, A., Pokomeda, K.:
	Application of response surface methodology and artificial neural network
	methods in modelling and optimization of biosorption process.
	\newblock Bioresource technology \textbf{160}, 150--160 (2014).
	\newblock \doi{10.1016/j.biortech.2014.01.021}
	
	\bibitem{Meireles.2003}
	Meireles, M., Almeida, P., Simoes, M.G.: A comprehensive review for industrial
	applicability of artificial neural networks.
	\newblock IEEE Transactions on Industrial Electronics \textbf{50}(3), 585--601
	(2003).
	\newblock \doi{10.1109/TIE.2003.812470}
	
	\bibitem{DelRioChanona.2017}
	{Del Rio-Chanona}, E.A., Fiorelli, F., Zhang, D., Ahmed, N.R., Jing, K., Shah,
	N.: An efficient model construction strategy to simulate microalgal lutein
	photo-production dynamic process.
	\newblock Biotechnology and Bioengineering \textbf{114}(11), 2518--2527 (2017).
	\newblock \doi{10.1002/bit.26373}
	
	\bibitem{Cheema.2002}
	Cheema, J.J.S., Sankpal, N.V., Tambe, S.S., Kulkarni, B.D.: Genetic programming
	assisted stochastic optimization strategies for optimization of glucose to
	gluconic acid fermentation.
	\newblock Biotechnology progress \textbf{18}(6), 1356--1365 (2002).
	\newblock \doi{10.1021/bp015509s}
	
	\bibitem{Desai.2008}
	Desai, K.M., Survase, S.A., Saudagar, P.S., Lele, S.S., Singhal, R.S.:
	Comparison of artificial neural network and response surface methodology in
	fermentation media optimization: Case study of fermentative production of
	scleroglucan.
	\newblock Biochemical Engineering Journal \textbf{41}(3), 266--273 (2008).
	\newblock \doi{10.1016/j.bej.2008.05.009}
	
	\bibitem{Nagata.2003}
	Nagata, Y., Chu, K.H.: Optimization of a fermentation medium using neural
	networks and genetic algorithms.
	\newblock Biotechnology Letters \textbf{25}(21), 1837--1842 (2003).
	\newblock \doi{10.1023/A:1026225526558}
	
	\bibitem{Fahmi.2012}
	Fahmi, I., Cremaschi, S.: Process synthesis of biodiesel production plant using
	artificial neural networks as the surrogate models.
	\newblock Computers {\&} Chemical Engineering \textbf{46}, 105--123 (2012).
	\newblock \doi{10.1016/j.compchemeng.2012.06.006}
	
	\bibitem{Nascimento.2000}
	Nascimento, C.A.O., Giudici, R., Guardani, R.: Neural network based approach
	for optimization of industrial chemical processes.
	\newblock Computers {\&} Chemical Engineering \textbf{24}(9-10), 2303--2314
	(2000).
	\newblock \doi{10.1016/S0098-1354(00)00587-1}
	
	\bibitem{Nascimento.1998}
	Nascimento, C.A.O., Giudici, R.: Neural network based approach for optimisation
	applied to an industrial nylon-6,6 polymerisation process.
	\newblock Computers {\&} Chemical Engineering \textbf{22}, 595--S600 (1998).
	\newblock \doi{10.1016/S0098-1354(98)00105-7}
	
	\bibitem{Chambers.2002}
	Chambers, M., Mount-Campbell, C.A.: Process optimization via neural network
	metamodeling.
	\newblock International Journal of Production Economics \textbf{79}(2), 93--100
	(2002).
	\newblock \doi{10.1016/S0925-5273(00)00188-2}
	
	\bibitem{Henao.2010}
	Henao, C.A., Maravelias, C.T.: Surrogate-based process synthesis.
	\newblock In: S.~Pierucci, G.B. Ferraris (eds.) 20th European Symposium on
	Computer Aided Process Engineering, \emph{Computer Aided Chemical
		Engineering}, vol.~28, pp. 1129--1134. Elsevier, Milano, Italy (2010).
	\newblock \doi{10.1016/S1570-7946(10)28189-0}
	
	\bibitem{Henao.2011}
	Henao, C.A., Maravelias, C.T.: Surrogate-based superstructure optimization
	framework.
	\newblock AIChE Journal \textbf{57}(5), 1216--1232 (2011).
	\newblock \doi{10.1002/aic.12341}
	
	\bibitem{SantAnna.2017}
	{Sant Anna}, H.R., Barreto, A.G., Tavares, F.W., de~Souza, M.B.: Machine
	learning model and optimization of a psa unit for methane-nitrogen
	separation.
	\newblock Computers {\&} Chemical Engineering \textbf{104}, 377--391 (2017).
	\newblock \doi{10.1016/j.compchemeng.2017.05.006}
	
	\bibitem{Smith.2013}
	Smith, J.D., Neto, A.A., Cremaschi, S., Crunkleton, D.W.: {CFD}-based
	optimization of a flooded bed algae bioreactor.
	\newblock Industrial {\&} Engineering Chemistry Research \textbf{52}(22),
	7181--7188 (2013).
	\newblock \doi{10.1021/ie302478d}
	
	\bibitem{Henao.2012}
	Henao, C.A.: A superstructure modeling framework for process synthesis using
	surrogate models.
	\newblock Dissertation, {University of Wisconsin}, Madison (2012)
	
	\bibitem{Kajero.2017}
	Kajero, O.T., Chen, T., Yao, Y., Chuang, Y.C., Wong, D.S.H.: Meta-modelling in
	chemical process system engineering.
	\newblock Journal of the Taiwan Institute of Chemical Engineers \textbf{73},
	135--145 (2017).
	\newblock \doi{10.1016/j.jtice.2016.10.042}
	
	\bibitem{Lewandowski.1998}
	Lewandowski, J., Lemcoff, N.O., Palosaari, S.: Use of neural networks in the
	simulation and optimization of pressure swing adsorption processes.
	\newblock Chemical Engineering {\&} Technology \textbf{21}(7), 593--597 (1998).
	\newblock \doi{10.1002/(SICI)1521-4125(199807)21:7<593::AID-CEAT593>3.0.CO;2-U}
	
	\bibitem{GutierrezAntonio.2016}
	Gutiérrez-Antonio, C.: Multiobjective stochastic optimization of dividing-wall
	distillation columns using a surrogate model based on neural networks.
	\newblock Chemical and Biochemical Engineering Quarterly \textbf{29}(4),
	491--504 (2016).
	\newblock \doi{10.15255/CABEQ.2014.2132}
	
	\bibitem{Chen.2002}
	Chen, C.R., Ramaswamy, H.S.: Modeling and optimization of variable retort
	temperature (vrt) thermal processing using coupled neural networks and
	genetic algorithms.
	\newblock Journal of Food Engineering \textbf{53}(3), 209--220 (2002).
	\newblock \doi{10.1016/S0260-8774(01)00159-5}
	
	\bibitem{Dornier.1995}
	Dornier, M., Decloux, M., Trystram, G., Lebert, A.: Interest of neural networks
	for the optimization of the crossflow filtration process.
	\newblock LWT-Food Science and Technology \textbf{28}(3), 300--309 (1995).
	\newblock \doi{10.1016/S0023-6438(95)94364-1}
	
	\bibitem{Fernandes.2006}
	Fernandes, F.A.N.: Optimization of fischer-tropsch synthesis using neural
	networks.
	\newblock Chemical Engineering {\&} Technology \textbf{29}(4), 449--453 (2006).
	\newblock \doi{10.1002/ceat.200500310}
	
	\bibitem{Grossmann.2002}
	Grossmann, I.E., Viswanathan, J., Vecchietti, A., Raman, R., Kalvelagen, E.:
	{GAMS/DICOPT}: A discrete continuous optimization package.
	\newblock GAMS Corporation Inc  (2002)
	
	\bibitem{Drud.1994}
	Drud, A.S.: Conopt - a large-scale grg code.
	\newblock ORSA Journal on Computing \textbf{6}(2), 207--216 (1994).
	\newblock \doi{10.1287/ijoc.6.2.207}
	
	\bibitem{Nandi.2001}
	Nandi, S., Ghosh, S., Tambe, S.S., Kulkarni, B.D.: Artificial
	neural-network-assisted stochastic process optimization strategies.
	\newblock AIChE Journal \textbf{47}(1), 126--141 (2001).
	\newblock \doi{10.1002/aic.690470113}
	
	\bibitem{Ryoo.1996}
	Ryoo, H.S., Sahinidis, N.V.: A branch-and-reduce approach to global
	optimization.
	\newblock Journal of Global Optimization \textbf{8}(2), 107--138 (1996).
	\newblock \doi{10.1007/BF00138689}
	
	\bibitem{Tawarmalani.2004}
	Tawarmalani, M., Sahinidis, N.V.: Global optimization of mixed-integer
	nonlinear programs: A theoretical and computational study.
	\newblock Mathematical Programming \textbf{99}(3), 563--591 (2004).
	\newblock \doi{10.1007/s10107-003-0467-6}
	
	\bibitem{Tawarmalani.2005}
	Tawarmalani, M., Sahinidis, N.V.: A polyhedral branch-and-cut approach to
	global optimization.
	\newblock Mathematical Programming \textbf{103}(2), 225--249 (2005).
	\newblock \doi{10.1007/s10107-005-0581-8}
	
	\bibitem{Bradford.2017}
	Bradford, E., Schweidtmann, A.M., Lapkin, A.A.: Efficient multiobjective
	optimization employing gaussian processes, spectral sampling and a genetic
	algorithm.
	\newblock In Revision  (2018)
	
	\bibitem{Forrester.2009}
	Forrester, A.I., Keane, A.J.: Recent advances in surrogate-based optimization.
	\newblock Progress in Aerospace Sciences \textbf{45}(1-3), 50--79 (2009).
	\newblock \doi{10.1016/j.paerosci.2008.11.001}
	
	\bibitem{Jones.2001}
	Jones, D.R.: A taxonomy of global optimization methods based on response
	surfaces.
	\newblock Journal of Global Optimization \textbf{21}(4), 345--383 (2001).
	\newblock \doi{10.1023/A:1012771025575}
	
	\bibitem{Sacks.1989}
	Sacks, J., Welch, W.J., Mitchell, T.J., Wynn, H.P.: Design and analysis of
	computer experiments.
	\newblock Statistical science \textbf{4}(4), 409--423 (1989).
	\newblock \doi{10.1214/ss/1177012413}
	
	\bibitem{Shahriari.2016}
	Shahriari, B., Swersky, K., Wang, Z., Adams, R.P., de~Freitas, N.: Taking the
	human out of the loop: A review of bayesian optimization.
	\newblock Proceedings of the IEEE \textbf{104}(1), 148--175 (2016).
	\newblock \doi{10.1109/JPROC.2015.2494218}
	
	\bibitem{Schweidtmann.2018}
	Schweidtmann, A.M., Holmes, N., Bradford, E., Bourne, R., Lapkin, A.:
	Multi-objective self-optimization of continuous flow reactors.
	\newblock In preparation  (2018)
	
	\bibitem{Wilson.2017}
	Wilson, Z.T., Sahinidis, N.V.: The {ALAMO} approach to machine learning.
	\newblock Computers {\&} Chemical Engineering \textbf{106}, 785--795 (2017).
	\newblock \doi{10.1016/j.compchemeng.2017.02.010}
	
	\bibitem{Cozad.2015}
	Cozad, A., Sahinidis, N.V., Miller, D.C.: A combined first-principles and
	data-driven approach to model building.
	\newblock Computers {\&} Chemical Engineering \textbf{73}, 116--127 (2015).
	\newblock \doi{10.1016/j.compchemeng.2014.11.010}
	
	\bibitem{Cozad.2014}
	Cozad, A., Sahinidis, N.V., Miller, D.C.: Learning surrogate models for
	simulation-based optimization.
	\newblock AIChE Journal \textbf{60}(6), 2211--2227 (2014).
	\newblock \doi{10.1002/aic.14418}
	
	\bibitem{Falk.1969}
	Falk, J.E., Soland, R.M.: An algorithm for separable nonconvex programming
	problems.
	\newblock Management Science \textbf{15}(9), 550--569 (1969).
	\newblock \doi{10.1287/mnsc.15.9.550}
	
	\bibitem{Horst.1996}
	Horst, R., Tuy, H.: Global Optimization: Deterministic Approaches, 3 edn.
	\newblock Springer, Berlin, Heidelberg (1996).
	\newblock \doi{10.1007/978-3-662-03199-5}
	
	\bibitem{Misener.2014}
	Misener, R., Floudas, C.A.: {ANTIGONE}: Algorithms for continuous / integer
	global optimization of nonlinear equations.
	\newblock Journal of Global Optimization \textbf{59}(2), 503--526 (2014).
	\newblock \doi{10.1007/s10898-014-0166-2}
	
	\bibitem{Maher.2017}
	Maher, S.J., Fischer, T., Gally, T., Gamrath, G., Gleixner, A., Gottwald, R.L.,
	Hendel, G., Koch, T., L{\"u}bbecke, M.E., Miltenberger, M., Müller, B.,
	Pfetsch, M.E., Puchert, C., Rehfeldt, D., Schenker, S., Schwarz, R., Serrano,
	F., Shinano, Y., Weninger, D., Witt, J.T., Witzig, J.: The {SCIP}
	optimization suite (version 4.0)
	
	\bibitem{Epperly.1997}
	Epperly, T.G.W., Pistikopoulos, E.N.: A reduced space branch and bound
	algorithm for global optimization.
	\newblock Journal of Global Optimization \textbf{11}(3), 287--311 (1997).
	\newblock \doi{10.1023/A:1008212418949}
	
	\bibitem{Mitsos.2009}
	Mitsos, A., Chachuat, B., Barton, P.I.: {McC}ormick-based relaxations of
	algorithms.
	\newblock SIAM Journal on Optimization \textbf{20}(2), 573--601 (2009).
	\newblock \doi{10.1137/080717341}
	
	\bibitem{Scott.2011}
	Scott, J.K., Stuber, M.D., Barton, P.I.: Generalized {McC}ormick relaxations.
	\newblock Journal of Global Optimization \textbf{51}(4), 569--606 (2011).
	\newblock \doi{10.1007/s10898-011-9664-7}
	
	\bibitem{Bongartz.2017}
	Bongartz, D., Mitsos, A.: Deterministic global optimization of process
	flowsheets in a reduced space using {McCormick} relaxations.
	\newblock Journal of Global Optimization \textbf{20}(9), 419 (2017).
	\newblock \doi{10.1007/s10898-017-0547-4}
	
	\bibitem{Huster.2017}
	Huster, W.R., Bongartz, D., Mitsos, A.: Deterministic global optimization of
	the design of a geothermal organic rankine cycle.
	\newblock Energy Procedia \textbf{129}, 50--57 (2017).
	\newblock \doi{10.1016/j.egypro.2017.09.181}
	
	\bibitem{McCormick.1976}
	McCormick, G.P.: Computability of global solutions to factorable nonconvex
	programs: Part {I} - convex underestimating problems.
	\newblock Mathematical Programming \textbf{10}(1), 147--175 (1976).
	\newblock \doi{10.1007/BF01580665}
	
	\bibitem{McCormick.1983}
	McCormick, G.P.: Nonlinear programming: theory, algorithms, and applications.
	\newblock {John Wiley {\&} Sons, Inc} (1983)
	
	\bibitem{Bompadre.2012}
	Bompadre, A., Mitsos, A.: Convergence rate of {McC}ormick relaxations.
	\newblock Journal of Global Optimization \textbf{52}(1), 1--28 (2012).
	\newblock \doi{10.1007/s10898-011-9685-2}
	
	\bibitem{Najman.2016}
	Najman, J., Mitsos, A.: Convergence analysis of multivariate
	{McCormick}relaxations.
	\newblock Journal of Global Optimization \textbf{66}(4), 597--628 (2016).
	\newblock \doi{10.1007/s10898-016-0408-6}
	
	\bibitem{Tsoukalas.2014}
	Tsoukalas, A., Mitsos, A.: Multivariate {McC}ormick relaxations.
	\newblock Journal of Global Optimization \textbf{59}(2-3), 633--662 (2014).
	\newblock \doi{10.1007/s10898-014-0176-0}
	
	\bibitem{Wechsung.2014}
	Wechsung, A., Barton, P.I.: Global optimization of bounded factorable functions
	with discontinuities.
	\newblock Journal of Global Optimization \textbf{58}(1), 1--30 (2014).
	\newblock \doi{10.1007/s10898-013-0060-3}
	
	\bibitem{Khan.2017}
	Khan, K.A., Watson, H.A.J., Barton, P.I.: Differentiable {McC}ormick
	relaxations.
	\newblock Journal of Global Optimization \textbf{67}(4), 687--729 (2017).
	\newblock \doi{10.1007/s10898-016-0440-6}
	
	\bibitem{Bongartz.2017b}
	Bongartz, D., Mitsos, A.: Infeasible path global flowsheet optimization using
	{McCormick} relaxations.
	\newblock In: A.~Espuna (ed.) 27th European Symposium on Computer Aided Process
	Engineering, \emph{Computer Aided Chemical Engineering}, vol.~40. Elsevier,
	San Diego (2017).
	\newblock \doi{10.1016/B978-0-444-63965-3.50107-0}
	
	\bibitem{Wechsung.2015}
	Wechsung, A., Scott, J.K., Watson, H.A.J., Barton, P.I.: Reverse propagation of
	{McC}ormick relaxations.
	\newblock Journal of Global Optimization \textbf{63}(1), 1--36 (2015).
	\newblock \doi{10.1007/s10898-015-0303-6}
	
	\bibitem{Stuber.2015}
	Stuber, M.D., Scott, J.K., Barton, P.I.: Convex and concave relaxations of
	implicit functions.
	\newblock Optimization Methods and Software \textbf{30}(3), 424--460 (2015).
	\newblock \doi{10.1080/10556788.2014.924514}
	
	\bibitem{Smith.1997}
	Smith, E.M., Pantelides, C.C.: Global optimisation of nonconvex minlps.
	\newblock Computers {\&} Chemical Engineering \textbf{21}, 791--796 (1997).
	\newblock \doi{10.1016/S0098-1354(97)87599-0}
	
	\bibitem{Tawarmalani.2002}
	Tawarmalani, M., Sahinidis, N.V., Pardalos, P.: Convexification and Global
	Optimization in Continuous and Mixed-Integer Nonlinear Programming: Theory,
	Algorithms, Software, and Applications, \emph{Nonconvex Optimization and Its
		Applications}, vol.~65.
	\newblock Springer, Boston, MA (2002).
	\newblock \doi{10.1007/978-1-4757-3532-1}
	
	\bibitem{Bishop.2009}
	Bishop, C.M.: Pattern recognition and machine learning, 8 edn.
	\newblock Information science and statistics. Springer, New York, NY (2009)
	
	\bibitem{Moore.1979}
	Moore, R.E., Bierbaum, F.: Methods and applications of interval analysis, 2
	edn.
	\newblock SIAM studies in applied mathematics. {Soc. for Industrial and Applied
		Mathematics}, Philadelphia (1979).
	\newblock \doi{10.1137/1.9781611970906}
	
	\bibitem{Hofschuster.1998}
	Hofschuster, W., Krämer, W.: {FILIB}++ interval library (version 3.0.2) (1998)
	
	\bibitem{Chachuat.2014}
	Chachuat, B.: {MC}++ (version 2.0): A toolkit for bounding factorable functions
	(2014)
	
	\bibitem{Chachuat.2015}
	Chachuat, B., Houska, B., Paulen, R., Peri'c, N., Rajyaguru, J., Villanueva,
	M.E.: Set-theoretic approaches in analysis, estimation and control of
	nonlinear systems.
	\newblock IFAC-PapersOnLine \textbf{48}(8), 981--995 (2015).
	\newblock \doi{10.1016/j.ifacol.2015.09.097}
	
	\bibitem{Bertsekas.2003}
	Bertsekas, D.P., Nedic, A., Ozdaglar, A.E.: Convex analysis and optimization,
	\emph{Athena Scientific optimization and computation series}, vol.~1.
	\newblock {Athena Scientific}, Belmont, Mass. (2003)
	
	\bibitem{Bertsekas.2015}
	Bertsekas, D.P.: Convex optimization algorithms, \emph{Optimization and
		computation series}, vol.~4.
	\newblock {Athena Scientific}, Belmont, Massachusetts (2015)
	
	\bibitem{Bongartz.2018}
	Bongartz, D., Najman, J., Scholl, S., Mitsos, A.: {MAiNGO}: {M}c{C}ormick based
	{A}lgorithm for mixed {i}nteger {N}onlinear {G}lobal {O}ptimization.
	\newblock Technical report  (2018)
	
	\bibitem{InternationalBusinessMachies.2009}
	{International Business Machies}: {IBM} ilog {CPLEX} (version 12.1) (2009)
	
	\bibitem{Gleixner.2017}
	Gleixner, A.M., Berthold, T., Müller, B., Weltge, S.: Three enhancements for
	optimization-based bound tightening.
	\newblock Journal of Global Optimization \textbf{67}(4), 731--757 (2017).
	\newblock \doi{10.1007/s10898-016-0450-4}
	
	\bibitem{Ryoo.1995}
	Ryoo, H.S., Sahinidis, N.V.: Global optimization of nonconvex {NLP}s and
	{MINLP}s with applications in process design.
	\newblock Computers {\&} Chemical Engineering \textbf{19}(5), 551--566 (1995).
	\newblock \doi{10.1016/0098-1354(94)00097-2}
	
	\bibitem{Locatelli.2013}
	Locatelli, M., Schoen, F. (eds.): Global optimization: Theory, algorithms, and
	applications.
	\newblock MOS-SIAM series on optimization. {Mathematical Programming Society},
	Philadelphia, Pa. (2013)
	
	\bibitem{Kraft.1988}
	Kraft, D.: A Software Package for Sequential Quadratic Programming.
	\newblock Deutsche Forschungs- und Versuchsanstalt für Luft- und Raumfahrt
	K{\"o}ln: Forschungsbericht. {Wiss. Berichtswesen d. DFVLR}, K{\"o}ln (1988)
	
	\bibitem{Kraft.1994}
	Kraft, D.: Algorithm 733: Tomp-fortran modules for optimal control
	calculations.
	\newblock ACM Transactions on Mathematical Software \textbf{20}(3), 262--281
	(1994).
	\newblock \doi{10.1145/192115.192124}
	
	\bibitem{Johnson.2016}
	Johnson, S.G.: The {NLopt} nonlinear-optimization package (version 2.4.2)
	(2016)
	
	\bibitem{Bendtsen.2012}
	Bendtsen, C., Stauning, O.: Fadbad++ (version 2.1): a flexible {C++} package
	for automatic differentiation (2012)
	
	\bibitem{Najman.2017b}
	Najman, J., Mitsos, A.: Tighter {McC}ormick relaxations through subgradient
	propagation.
	\newblock In Revision \url{https://arxiv.org/abs/1710.09188} (2017)
	
	\bibitem{Ghorbanian.2009}
	Ghorbanian, K., Gholamrezaei, M.: An artificial neural network approach to
	compressor performance prediction.
	\newblock Applied Energy \textbf{86}(7-8), 1210--1221 (2009).
	\newblock \doi{10.1016/j.apenergy.2008.06.006}
	
	\bibitem{Luyben.2010}
	Luyben, W.L.: Design and control of the cumene process.
	\newblock Industrial {\&} Engineering Chemistry Research \textbf{49}(2),
	719--734 (2010).
	\newblock \doi{10.1021/ie9011535}
	
	\bibitem{Schultz.2016}
	Schultz, E.S., Trierweiler, J.O., Farenzena, M.: The importance of nominal
	operating point selection in self-optimizing control.
	\newblock Industrial {\&} Engineering Chemistry Research \textbf{55}(27),
	7381--7393 (2016).
	\newblock \doi{10.1021/acs.iecr.5b02044}
	
	\bibitem{Lee.2016}
	Lee, U., Burre, J., Caspari, A., Kleinekorte, J., Schweidtmann, A.M., Mitsos,
	A.: Techno-economic optimization of a green-field post-combustion {CO} 2
	capture process using superstructure and rate-based models.
	\newblock Industrial {\&} Engineering Chemistry Research \textbf{55}(46),
	12,014--12,026 (2016).
	\newblock \doi{10.1021/acs.iecr.6b01668}
	
	\bibitem{Helmdach.2017}
	Helmdach, D., Yaseneva, P., Heer, P.K., Schweidtmann, A.M., Lapkin, A.A.: A
	multiobjective optimization including results of life cycle assessment in
	developing biorenewables-based processes.
	\newblock ChemSusChem \textbf{10}(18), 3632--3643 (2017).
	\newblock \doi{10.1002/cssc.201700927}
	
\end{thebibliography}

\end{document}